\documentclass[a4paper,12pt]{amsart}
\usepackage[english,french]{babel}
\usepackage{amsmath, amsfonts, amssymb, amsthm,mathptmx}
\usepackage{color}
\usepackage{graphicx}
\usepackage{multirow}

\usepackage{rotating,pdflscape,color,amssymb,hyperref,subfigure,psfrag,amsmath,eufrak,bbm,epsfig}

\newtheorem{lemme}{Lemma}[section]

\newtheorem{theo}[lemme]{Theorem}

\newenvironment{dem*}{\underline{\textbf{Proof }}}{\hfill$\square$ $$$$}

\newtheorem{assum}[lemme]{Assumption} 
\theoremstyle{definition}

\newtheorem{rem}[lemme]{Remark}

\DeclareMathOperator{\im}{\operatorname{Im}}
\DeclareMathOperator{\re}{\operatorname{Re}}

\DeclareMathOperator{\eloi}{\overset{(d)}{=}}

\DeclareMathOperator{\cloi}{\overset{(d)}{\longrightarrow}}
\DeclareMathOperator{\cvp}{\overset{(\mathbb{P})}{\longrightarrow}}
\DeclareMathOperator{\E}{\mathbb{E}}

\DeclareMathOperator{\pro}{\mathbb{P}}

\DeclareMathOperator{\op}{op}

\DeclareMathOperator{\diag}{\text{diag}}

\DeclareMathOperator{\M}{\mathcal{M}}
\DeclareMathOperator{\C}{\mathbb{C}}

\DeclareMathOperator{\N}{\mathbb{N}}

\newcommand{\ope}{\operatorname}
\newcommand{\tto}{\longrightarrow}
\newcommand{\ol}{\overline}
\newcommand{\wt}{\widetilde}

\newcommand{\ep}{\varepsilon}

\newcommand{\wg}{\operatorname{Wg}}

\newcommand{\card}{\operatorname{Card}}

\newcommand{\tr}{\operatorname{Tr}}
\newcommand{\dI}{\mathrm{I}}
\newcommand{\dV}{\mathrm{V}}
\newcommand{\bA}{\mathbf{A}}

\newcommand{\bU}{\mathbf{U}}
\newcommand{\bT}{\mathbf{T}}

\newcommand{\bB}{\mathbf{B}}
\newcommand{\bC}{\mathbf{C}}
\newcommand{\bD}{\mathbf{D}}
\newcommand{\bI}{\mathbf{I}}

\newcommand{\bQ}{\mathbf{Q}}
\newcommand{\bR}{\mathbf{R}}
\newcommand{\bv}{\mathbf{v}}
\newcommand{\be}{\mathbf{e}}

\newcommand{\bu}{\mathbf{u}}
\newcommand{\bM}{\mathbf{M}}

\newcommand{\bW}{\mathbf{W}}

\newcommand{\bH}{\mathbf{H}}

\newcommand{\bJ}{\mathbf{J}}
\newcommand{\bY}{\mathbf{Y}}

\newcommand{\bZ}{\mathbf{Z}}

\newcommand{\tred}{\textcolor{red}}

\newcommand{\bpr}{\begin{pr}}
\newcommand{\epr}{\end{pr}}
\newcommand{\trm}{\textrm}
\newcommand{\f}{\frac}

\newcommand{\ff}{\frac{1}}
\newcommand{\one}{\mathbbm{1}}

\newcommand{\ds}{\displaystyle}
\newcommand{\lf}{\left}
\newcommand{\ri}{\right}
\newcommand{\st}{such that }

\newcommand{\dist}{\operatorname{dist}}
\newcommand{\bes}{\begin{equation*}}
\newcommand{\ees}{\end{equation*}}
\newcommand{\beqy}{\begin{eqnarray}}
\newcommand{\eeqy}{\end{eqnarray}}
\newcommand{\beq}{\begin{eqnarray*}}
\newcommand{\eeq}{\end{eqnarray*}}
\newcommand{\bbe}{\begin{equation}}
\newcommand{\ee}{\end{equation}}
\newcommand{\bbm}{\begin{bmatrix}}
\newcommand{\ebm}{\end{bmatrix}}
\newcommand{\bpm}{\begin{pmatrix}}
\newcommand{\epm}{\end{pmatrix}}
\newcommand{\bdet}{\begin{vmatrix}}
\newcommand{\edet}{\end{vmatrix}}
\newcommand{\la}{\label}
\newcommand{\eqre}{\eqref}
\newcommand{\ti}{\times}

\newcommand{\egd}{\ := \ }

\newcommand{\OO}[1]{O \left( #1 \right)}
\newcommand{\oo}[1]{o \left( #1 \right)}
 
\newcommand{\ovl}{\overline}
\newcommand{\alp}{\alpha}
\newcommand{\al}{\alp}

\newcommand{\bet}{\beta}

\newcommand{\tta}{\theta}
\newcommand{\lam}{\lambda}
\newcommand{\si}{\sigma}

\newcommand{\ld}{\ldots}
\newcommand{\norm}[1]{\left\| #1 \right\|_{\op}}

\newcommand{\lexp}[2]{{\vphantom{#2}}^{#1}#2}
\newcommand{\lbinom}[3]{{\vphantom{#3}}^{\;\;#1}_{#2}#3}

\newcommand{\Ec}[1]{\E \left[ #1 \right]}
\newcommand{\inve}[1]{ \left( #1 \right)^{-1}}

\newcommand{\bgt}{\begin{itemize}}
\newcommand{\ent}{\end{itemize}}
\newcommand{\ite}{\item}

\newenvironment{pr}{\noindent {\bf Proof. }}{\hfill $\square$\\}

\newcommand{\spec}{\operatorname{Spec}}
\newcommand{\scalv}[2]{\big<#1 , \ #2 \big>}
\newcommand{\ii}{\operatorname{i}}
\newcommand{\supp}{\operatorname{supp}}

\begin{document}
\title{Complex outliers of  Hermitian random matrices}
\author{Jean Rochet}


\thanks{JR: MAP5,
Universit\'e Paris Descartes,
45, rue des Saints-P\`eres
75270 Paris Cedex 06, France.  jean.rochet@parisdescartes.fr.}

\maketitle

 \renewcommand{\abstractname}{Abstract}

\begin{abstract}
In this paper, we study the asymptotic behavior of the outliers of the sum a Hermitian random matrix and a finite rank matrix which is not necessarily Hermitian. We observe several possible  convergence rates  and  outliers locating around their limits at the vertices of regular polygons  as in \cite{FloJean}, as well as possible correlations between outliers at macroscopic distance as in \cite{KYIOLW} and \cite{FloJean}. We also observe that a single spike can generate  several outliers in the spectrum of the deformed model, as already noticed in \cite{BEN1} and \cite{BEL}.  In the particular case where the perturbation matrix is Hermitian, our results complete the work of \cite{BEN2}, as we consider fluctuations of outliers lying in ``holes" of the limit support,   which happen to exhibit surprising correlations.\end{abstract}

\section{Introduction}
It is known that   adding a finite rank perturbation to a large matrix barely changes the global behavior of its spectrum. Nevertheless, some of the eigenvalues, called \emph{outliers}, can deviate away from the bulk of the spectrum, depending on the strength of the perturbation. This phenomenon, well known as the \emph{BBP transition}, was first brought to light for empirical covariance matrices by Johnstone in \cite{john:01}, by Baik, Ben Arous and P\'ech\'e in \cite{BAI}, and then shown under several hypothesis in the Hermitian case in \cite{SP06,FP07,CAP2,CAP,PRS11, BEN1,BEN4,BEN2,BEN3,CAP3,KYIsoSCL, KYIOLW}. Non-Hermitian models have been also studied: i.i.d. matrices in \cite{TAO,BORCAP1,RAJ}, elliptic matrices in \cite{ROR1} and matrices from the Single Ring Theorem in \cite{FloJean}. In \cite{FloJean}, and lately in \cite{RAJ}, the authors have also studied the fluctuations of the outliers and, due to non-Hermitian structure, obtained unusual results: the distribution of the fluctuations highly depends on the shape of the Jordan Canonical Form of the perturbation, in particular, the convergence rate depends on the size of the Jordan blocks. Also, the outliers tend to locate around their limit at the vertices of a regular polygon. At last, they observe correlations between the fluctuations of outliers at a macroscopic distance with each other. \\
In this paper, we show that the same kind of phenomenon occurs when we perturb an Hermitian matrix $\bH$ with a non-Hermitian one $\bA$. More precisely, we study finite rank perturbations for Hermitian random matrices $\bH$ whose spectral measure tends to a compactly supported measure $\mu$ and the perturbation $\bA$ is just a complex matrix with a finite rank. With further assumptions, we prove that outliers of $\bH+\bA$ may appear at a macroscopic distance from the bulk and, following the ideas of \cite{FloJean}, we show that they fluctuate with convergence rates which depend on the matrix $\bA$ through its Jordan Canonical Form. Remind that any   complex matrix is similar to a   block diagonal matrix with diagonal blocks  of the type
 $$\bR_p(\tta) \ := \ \bpm \tta & 1 & & (0) \\ 
                             &\ddots&\ddots &\\ 
                             \multicolumn{2}{c}{\multirow{2}{*}{(0)}} &\ddots & 1 \\
                             \multicolumn{2}{c}{}                     & &\tta \epm \\,
                             $$
so that $\bA \sim \diag\big(\bR_{p_1}(\tta_1),\ld,\bR_{p_q}(\tta_q) \big)$, this last matrix being called the \emph{Jordan Canonical Form} of   $\bA$ \cite[Chapter 3]{HORNJOHNSON}. We show, up to some hypothesis, that for any eigenvalue $\tta$ of $\bA$, if we denote by
$$
\underbrace{p_1,\ld,p_1}_{\bet_1 \text{ times}} > \underbrace{p_2, \ld , p_2}_{\bet_2 \text{ times}} > \cdots > \underbrace{p_\al, \ld, p_\al}_{\bet_\al \text{ times}} 
$$
the sizes of the blocks associated to $\tta$   in the Jordan Canonical Form of $\bA$ and  introduce the (possibly empty) set
$$
\mathcal{S}_\tta \ : = \ \lf\{ \xi \in \C, \ G_\mu(\xi) = \ff \tta \ri\}
$$
where $G_\mu(z) := \ds\int \ff{z-x}\mu(dx)$ is the Cauchy transform of the measure $\mu$, then there are exactly $\bet_1 p_1 + \cdots + \bet_\al p_\al$ outliers of $\bH+\bA$ tending to each element of $\mathcal{S}_\tta$. We also prove that   for each element $\xi$ in $\mathcal{S}_\tta$, there are exactly $\bet_1 p_1$ outliers  tending to $\xi$ at rate $N^{-1/(2p_1)}$, $\bet_2 p_2$ outliers tending to $\xi$ at rate $N^{-1/(2p_2)}$, etc... (see Figure \ref{wignerfluc1348012015}).
Furthermore, the limit joint distribution of the fluctuations is explicit, not necessarily Gaussian, and might show  correlations even between outliers at a macroscopic distance with each other. 
This phenomenon of correlations between the fluctuations of two   outliers with distinct limits has already been proved for non-Gaussian Wigner matrices when $\bA$ is Hermitian (see \cite{KYIOLW}), while in our case, Gaussian Wigner matrices can have such correlated outliers: indeed, the correlations that we bring to light here are due to the fact that the eigenspaces of $\bA$ are not necessarily orthogonal or that one single spike generates several outliers. 
Indeed, we observe that the outliers may outnumber the rank of $\bA$. This had   already been noticed  in \cite[Remark 2.11]{BEN1} when the support of the limit spectral measure of $\bH$ has some ``holes'' or in the   different model of  \cite{BEL}, where the authors study the case where $\bA$ is Hermitian but with full rank and is invariant  in distribution by unitary conjugation. Here, the phenomenon can be proved to occur even when the   support of the limit spectral measure of $\bH$ is connected.
At last, if we apply our results in the particular case where $\bA$ is Hermitian, we also see that two outliers at a macroscopic distance with each other are correlated if they both are generated by the same spike (which can occur only if the limit support   is  disconnected) and are independent otherwise (see Figure \ref{refFigure17154567}). From this point of view, this completes the work of \cite{BEN2}, where fluctuations of outliers lying in ``holes" of the limit support had not been studied.\\
The fact to consider a non-Hermitian deformation on a Hermitian random matrix has already been studied in theoretical physics (see \cite{FYO1,FYO2,FYO3,FYO4}) in the particular case where $\bH$ is a GOE/GUE matrix and $\bA$ is a non negative Hermitian matrix times $\ii$ (the square root of $-1$). They proved a weaker version of Theorem \ref{theo1123201504058899} in this specific case but didn't study the fluctuations. \\
The proofs of this paper rely essentially on the ideas of the paper \cite{FloJean} about outliers in the Single Ring Theorem and on the results proved in  \cite{PRS11,PRS12,BEN2}. More precisely, the study of the fluctuations reproduce the outlines of the proofs of \cite{FloJean} as long as the model fulfills some conditions. Thanks to \cite{PRS11,PRS12}, we show that these conditions are satisfied for Wigner matrices. At last, using \cite{BEN2} and the Weingarten calculus, we show the same for Hermitian matrices invariant in distribution by unitary conjugation. In the appendix, as a tool for the outliers study, we prove a result on the fluctuations of the entries of such matrices.

%


%

\section{General Results}

At first, we formulate the results in general settings and we shall give, in the next section, examples of random matrices on which these results apply.

\subsection{Convergence of the outliers}

\subsubsection{Set up and assumptions} For all $N \geq 1$, let $\bH_N$ be an Hermitian random $N \ti N$ matrix whose \emph{empirical spectral measure}, as $N$ goes to infinity, converges weakly in probability to a compactly supported measure $\mu$
\beqy \la{assu1055m2015k}
\mu_N & : = & \ff N \sum_{i=1}^N \delta_{\lam_i(\bH)} \ \tto \ \mu.
\eeqy
We shall suppose that $\mu$ is non trivial in the sense that $\mu$ is not a single Dirac measure. Also, we suppose that $\bH_N$ does not possess any \emph{natural outliers}, i.e.
 \begin{assum}\la{assu180203062015}
 As $N$ goes to infinity, with probability tending to one, 
\beq
\sup_{\lam \in \spec(\bH_N)} \dist(\lam,\supp(\mu)) & \tto & 0.
\eeq
\end{assum}
For all $N \geq 1$, let $\bA_N$ be an $N \ti N$ random matrix independent from $\bH_N$ (which does not satisfies necessarily $\bA_N^* = \bA_N$) whose  rank is bounded by an integer $r$ (independent from $N$).  We know that we can write \beqy \label{105812052015assu}
\bA_N & := & \bU \bpm \bA_0 & 0 \\ 0 & 0 \\ \epm \bU^* 
\eeqy 
where $\bU$ is an $N\ti N$ unitary matrix and $\bA_0$ is $2r \ti 2r$ matrix. We notice that $\bA_N$ only depends on the $2r$-first columns of $\bU$ so that, we shall write 
$$
\bA_N \ := \ \bU_{2r} \bA_0 \bU_{2r}^*,
$$ 
where the $N\ti 2r$ matrix $\bU_{2r}$ designates the $2r$-first columns of $\bU$. We shall assume that $\bA_0$ is deterministic and independent from $N$. We shall denote by $\tta_1,\ld,\tta_j$ the distinct non-zero eigenvalues of $\bA_0$ and $k_1,\ld,k_j$ their respective multiplicity\footnote{The \emph{multiplicity} of an eigenvalue is defined as its order as a root of the characteristic polynomial, which is greater than or equal to the dimension of the associated eigenspace.} (note that $\sum_{i=1}^jk_i \leq r$).
We consider the additive perturbation
\beqy
\wt\bH_N & := &  \bH_N + \bA_N,
\eeqy

\noindent We set
\beqy
G_\mu(z)& := & \int \ff{z-x}\mu(dx) .
\eeqy 
the Cauchy transform of the measure $\mu$. We introduce, for all $i \in \{1,\ld,j\}$, the finite, possibly empty, set
\beqy\la{setdef122720150206}
\mathcal{S}_{\tta_i} &  : = & \lf\{ \xi \in \C\backslash \supp(\mu), \ G_\mu(\xi) = \ff \tta_i \ri\}, \ \text{ and } \ m_i \ : = \ \card \mathcal{S}_{\tta_i}
\eeqy


We make the following assumption

\begin{assum}\la{assu173602062015}
For any $\delta>0$, as $N$ goes to infinity, we have
$$
\sup_{\dist(z,\supp(\mu))>\delta }\norm{\bU_{2r}^* \inve{z\bI - \bH_N} \bU_{2r} - G_\mu(z) \bI} \ \cvp \ 0.
$$
\end{assum}

\subsubsection{Result}

\begin{theo}[Convergence of the outliers] \la{theo1123201504058899}
For  $\tta_1,\ld,\tta_j$, $k_1,\ld, k_j$, $\mathcal{S}_{\tta_1},\ld,\mathcal{S}_{\tta_j}$ and $m_1,\ld,m_j$ as defined above, with probability tending to one, $\wt\bH_N := \bH_N + \bA_N$ possesses exactly $\ds\sum_{i=1}^j k_i m_i$ eigenvalues at a macroscopic distance of $\supp \mu$ (outliers). More precisely, for all small enough $\delta>0$, for all large enough $N$, for all $i \in \{1,\ld,j\}$, if we set 
$$
\mathcal{S}_{\tta_i} \ = \ \{\xi_{i,1},\ld,\xi_{i,m_i}\},
$$
there are $m_i$ eigenvalues $\wt\lam_{i,1},\ld,\wt\lam_{i,m_i}$ of $\wt\bH_N$ in $\{z, \ \dist(z,\supp(\mu))>\delta\}$ satisfying
\beq
\wt\lam_{i,n} & = & \xi_{i,n} + \oo1, \ \ \ \text{ for all } n \in \{1,\ld,m_i\},
\eeq
after a proper labeling.
\end{theo}

\begin{rem}
If all the $\mathcal{S}_{\tta_i}$'s are empty, there is possibly no outlier  at all. This condition is the analogous of the \emph{phase transition} condition in \cite[Theorem 2.1]{BEN1} in the case where the $\tta_i$'s are real, which is if
$$
 \ff \tta_i \not\in \lf] \lim_{x \to a^-}G_{\mu}(x) , \lim_{x \to b^+}G_{\mu}(x) \ri[
$$
where $a$ (resp. $b$) designates the infimum (resp. the supremum) of the support of $\mu$,  
then, $\tta_i$ does not generate any outlier. In our case, if $|\tta_i|$ is large enough, $\mathcal{S}_{\tta_i}$ is necessarily non-empty\footnote{due to the fact that the Cauchy transform of a compactly supported measure can always be inverted in a neighborhood of infinity.}, which means that a strong enough perturbation always creates outliers.
\end{rem}

\begin{figure}[ht]
\centering
\subfigure[Spectrum of an Hermitian matrix of size $N=2000$, whose spectral measure tends to the \emph{semi-circle law} $\mu_{sc}(dx) := \ff{2\pi}\sqrt{4-x^2}\one{[-2,2]}(dx)$ (such as Wigner matrix), with perturbation matrix $\bA = \diag(i\sqrt{2},0\ld,0)$. 
]{\includegraphics[scale=0.45]{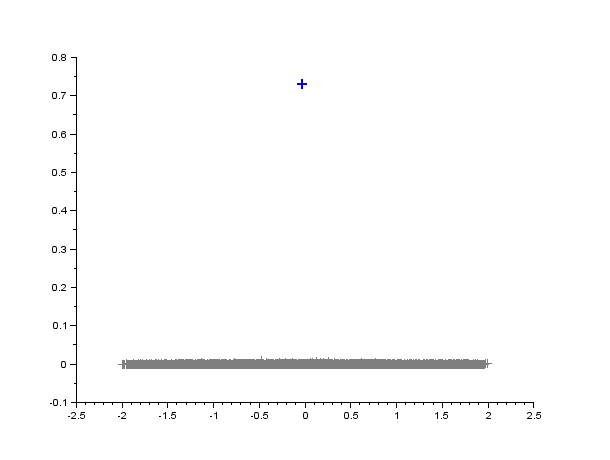}}\qquad
\subfigure[Spectrum of an Hermitian matrix of size $N=2000$, whose spectral measure tends to $\f{2}{5}\delta_{-1}(dx)+\f{2}{5}\delta_1(dx) + \ff{5}\mu_{sc}(dx)$ and with perturbation matrix $\bA = \diag(i\sqrt{2},0\ld,0)$.]{\includegraphics[scale=0.45]{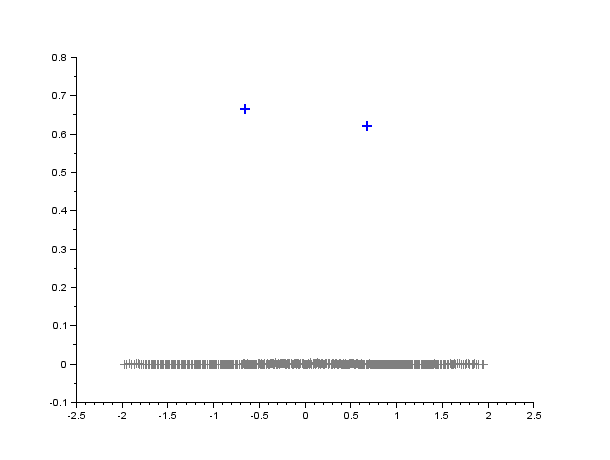} \label{contrexemple1435222015}}
\caption{Spectrums of two Hermitian matrices with the same limit bulk but different limit spectral densities on this bulk,   perturbed by the same matrix: both do not have the same number of outliers (the blue crosses ``\textcolor{blue}{$+$}'').}
\end{figure}

\begin{rem}
We notice that the outliers can outnumber the rank of $\bA$. This phenomenon was already observed in \cite{BEN1} in the case where the support of the limit spectral distribution has a disconnected support (see also \cite{BEL}). In our case, the phenomenon occurs even for connected support (see Figure \ref{contrexemple1435222015}).
\end{rem}

\subsection{Fluctuations of the outliers}
To study the fluctuations, one needs to understand the limit distribution of
 \beqy\la{truc1148254520150306}\sqrt N \norm{\bU_{2r}^* \inve{z\bI - \bH_N} \bU_{2r} - G_\mu(z) \bI}.\eeqy
 In the particular case where $\bH_N$ is a Wigner matrix, we know from \cite{PRS11} that this quantity is tight but does not necessarily converge. Hence, we shall need additional assumptions.

\subsubsection{Set up and assumptions}
As $\bA_N$ is not Hermitian, we need to introduce the Jordan Canonical Form (JCF) to describe the fluctuations. More precisely, we shall consider the JCF of $\bA_0$ which does not depend on $N$. We know that, in a proper basis, $\bA_0$ is a direct sum of \emph{Jordan blocks}, i.e. blocks of the form
\beqy
\bR_p(\tta) & = &  \bpm \tta & 1 & & (0) \\ 
                             &\ddots&\ddots &\\ 
                             \multicolumn{2}{c}{\multirow{2}{*}{(0)}} &\ddots & 1 \\
                             \multicolumn{2}{c}{}                     & &\tta \epm, \ p\ti p \text{ matrix, } \qquad \tta \in \C, p \geq 1
\eeqy

Let us denote by $\tta_1, \ld, \tta_q$ the distinct eigenvalues of $\bA_0 $ such that $\mathcal{S}_{\tta} \neq \emptyset$ (see \eqref{setdef122720150206} for the definition of $\mathcal{S}_\tta$), and for each $i=1, \ld, q$, we introduce a positive integer $\al_i$, some positive integers      $p_{i,1}> \cdots> p_{i,\al_i}$ corresponding to the distinct sizes of the blocks relative to the eigenvalue $\tta_i$  and $\bet_{i,1}, \ld, \bet_{i, \al_i}$ \st  for all $j$, $\bR_{p_{i,j}}(\tta_i)$ appears $\bet_{i,j}$ times, so that, for a certain $2r \ti 2r$ non singular matrix $\bQ$, we have:\\ 

\beqy\label{Eq14113004200158787} \bJ &=& \bQ^{-1}\bA_0\bQ  \ = \
  \hat\bA \ \bigoplus \ \bigoplus_{i=1}^q \ \bigoplus_{j=1}^{\alp_i} \! \underbrace{\begin{pmatrix}\bR_{p_{i,j}}(\tta_i)&&\\ &\ddots&\\ &&\bR_{p_{i,j}}(\tta_i)\end{pmatrix}}_{
 \bet_{{i,j}} \trm{ blocks}}
 \eeqy
where $\oplus$ is defined, for square block matrices, by $\mathbf{M}\oplus \mathbf{N}:=\bpm \mathbf{M}& 0\\ 0&\mathbf{N}\epm$ and $\hat\bA$ is a matrix such that its eigenvalues $\tta$ are such that $\mathcal{S}_\tta = \emptyset$ or null.\\ 
\indent The asymptotic orders of the fluctuations of the eigenvalues of $\bH_N+\bA_N $ depend  on the sizes $p_{i,j}$ of the blocks. Actually, for each $\tta_i$ and each $\xi_{i,n} \in \mathcal{S}_{\tta_i} = \{\xi_{i,1},\ld,\xi_{i,m_i}\}$,
 we know, by Theorem \ref{theo1123201504058899},  there are $\sum_{j=1}^{\alp_i}p_{ij}\ti\bet_{i,j}$ eigenvalues ${\wt\lambda}$ of $\bH_N+\bA_N$  which tend  to $\xi_{i,n}$ : we shall write them with a $\xi_{i,n}$ 
 on the top left corner, as follows $$\lexp{\xi_{i,n}}{{\wt\lambda}}.$$ Theorem \ref{belowtheor14122015} below will state that  for each block with size 
 $p_{i,j}$  corresponding to $\tta_i$  in  the JCF of $\bA_0$, there are $p_{i,j}$ eigenvalues (we shall write them with $p_{i,j}$ on the bottom left corner : $\lbinom{\xi_{i,n}}
 {p_{i,j}}{{\wt\lambda}}$) whose  convergence rate will be $N^{-1/(2p_{i,j})}$. As there are $\bet_{{i,j}}$ blocks of size 
 $p_{i,j}$, there are actually $p_{i,j}\times \bet_{{i,j}}$ eigenvalues tending to $\xi_{i,n}$ with  convergence rate    
 $N^{-1/(2p_{i,j})}$ (we shall write them $\lbinom{\xi_{i,n}}{p_{i,j}}{{\wt\lambda}_{s,t}}$ with $s \in \{1,\ldots,p_{i,j}\}$ and 
 $t \in \{1,\ldots,\bet_{{i,j}}\}$). It would be convenient to denote by $\Lambda_{i,j,n}$ the vector with size $p_{i,j}\times \bet_{{i,j}}$ defined by 
 \beqy \label{defLambda30042015}
 \Lambda_{i,j,n} \egd  \displaystyle \left(N^{1/(2p_{i,j})} \cdot \Big(\lbinom{\xi_{i,n}}{p_{i,j}}{ {\wt\lambda}_{s,t}} - \xi_{i,n} \Big) \right)_{\substack{1\le s\le p_{i,j}\\ 1\le t\le \bet_{i,j}}}.
\eeqy

In addition, we make an assumption on the convergence of \eqref{truc1148254520150306}.

\begin{assum}\la{lem15052015} \text{ }
\bgt
\ite[(1)]  The vector $\ds \lf(\sqrt{N}\bU_{2r}^*\lf( \inve{\xi_{i,n}-\bH_N} - \ff {\tta_i}\ri)\bU_{2r}\ri)_{\substack{1 \leq \; i\leq q_{\hphantom{1}} \\ 1 \leq n \leq m_i}}$ converges   in distribution and none of its entries tends to zero.
\ite[(2)] For all $k \geq 1$, all $i\in \{1,\ld,q\}$ and all $n \in \{1,\ld,m_i\}$, 
$$
\sqrt{N}\bU_{2r}^* \lf(\big(\xi_{i,n}-\bH_N\big)^{-(k+1)} - \int \f{\mu(dx)}{(\xi_{i,n}-x)^{k+1}} \ri)\bU_{2r}
$$
is tight.
\ent
or 
\bgt
\ite[(0')] For all $i\in \{1,\ld,q\}$ and all $n\in \{1,\ld,m_i\}$, as $N$ goes to infinity, $$\ds\sqrt{N}\lf(\ff N \tr \inve{\xi_{i,n}-\bH_N} - \ff{\tta_i}\big)\ri) \tto 0.$$
\ite[(1')] The vector $\ds \lf(\sqrt{N}\bU_{2r}^*\lf( \inve{\xi_{i,n}-\bH_N} - \ff N \tr\big(\xi_{i,n}-\bH_N\big)^{-1}\ri)\bU_{2r}\ri)_{\substack{1 \leq \; i\leq q_{\hphantom{1}} \\ 1 \leq n \leq m_i}}$ converges in distribution and none of its entries tends to zero.
\ite[(2')] For all $k \geq 1$ and for all $i\in \{1,\ld,q\}$, 
$$
\sqrt{N}\bU_{2r}^* \lf(\big(\xi_{i,n}-\bH_N\big)^{-(k+1)} -  \ff N \tr\big(\xi_{i,n}-\bH_N\big)^{-(k+1)}\ri)\bU_{2r}
$$
is tight.
\ent
\end{assum}

\indent As in \cite{FloJean}, we define now the family of random matrices that we shall use to characterize the limit distribution of the $\Lambda_{i,j,n}$'s. 
For each $i=1, \ldots, q$, let $I(\tta_i)$ (resp. $J(\tta_i)$) denote the set, with cardinality $\sum_{j=1}^{\al_i}\bet_{i,j}$, of indices in $\{1, \ld, r\}$ corresponding to the first (resp. last) columns of the blocks $\bR_{p_{i,j}}(\tta_i)$ ($1\le j\le \al_i$) in \eqre{Eq14113004200158787}.
\begin{rem} \la{rem1631290920142015}
 Note that the columns of $\bQ$ (resp. of $(\bQ^{-1})^*$) whose index belongs to $I(\tta_i)$ (resp. $J(\tta_i)$) are   eigenvectors of $\bA_0$  (resp. of $\bA_0^*$) associated to $\tta_i$ (resp. $\ol{\tta_i}$). See \cite[Remark 2.7]{FloJean}.
 \end{rem}

Now, let  \bbe\la{2121414h22015}
\lf({m}^{\tta_i,n}_{k,\ell}\ri)_{\substack{1 \leq i \leq q_{\hphantom{1.}} \quad\quad \ \ \\ 1 \leq n \leq m_i \quad\quad \ \ \\ (k,\ell)\in J(\tta_i)\ti I(\tta_i)}} \ee 
 be the multivariate random variable defined as the limit joint-distribution of
    \bbe\la{2121414h32015}
 \lf(\sqrt{N}\be_k^* \bQ^{-1}\bU_{2r}^*\lf( \inve{\xi_{i,n}-\bH_N} - \ff {\tta_i}\ri)\bU_{2r}\bQ\be_\ell\ri)_{\substack{1 \leq i \leq q_{\hphantom{1.}} \quad\quad \ \ \\ 1 \leq n \leq m_i \quad\quad \ \ \\ (k,\ell)\in J(\tta_i)\ti I(\tta_i)}} \ \underset{\text{jointly}}{\cloi} \ \lf({m}^{\tta_i,n}_{k,\ell}\ri)_{\substack{1 \leq i \leq q_{\hphantom{1.}} \quad\quad \ \ \\ 1 \leq n \leq m_i \quad\quad \ \ \\ (k,\ell)\in J(\tta_i)\ti I(\tta_i)}}
    \ee 
(which does exist by Assumption \ref{lem15052015}) and    where    $\be_1, \ld, \be_r$ are  the column vectors of the canonical basis of $\C^r$).

 For each $i,j$, let $K(i,j)$ (resp. $K(i,j)^-$) be the set, with cardinality $\bet_{i,j}$ (resp. $\sum_{j'=1}^{j-1}\bet_{i,j'}$), of  indices in $J(\tta_i)$    corresponding to a block of the type $\bR_{p_{i,j}}(\tta_i)$ (resp. to a block of the type  $\bR_{p_{i,j'}}(\tta_i)$ for  $j'<j$). In the same way, let $L(i,j)$ (resp. $L(i,j)^-$) be the set, with the same cardinality as $K(i,j)$ (resp. as $K(i,j)^-$), of indices in $I(\tta_i)$ corresponding to a block of the type $\bR_{p_{i,j}}(\tta_i)$ (resp. to a block of the type  $\bR_{p_{i,j'}}(\tta_i)$ for  $j'<j$).  Note that $K(i,j)^-$ and $L(i,j)^-$ are empty if $j=1$.  Let us define the random matrices for each $n \in \{1,\ld,m_i\}$
 \beqy \label{section2.32015}
 \ope{M}^{\tta_i,\mathrm{I}}_{j,n}\egd[m^{\tta_i,n}_{k,\ell}]_{\ds^{k\in K(i,j)^-}_{\ell \in L(i,j)^-}} &\qquad \qquad \qquad& \ope{M}^{\tta_i,\dI\dI}_{j,n}\egd[m^{\tta_i,n}_{k,\ell}]_{\ds^{k\in K(i,j)^-}_{\ell\in L(i,j)}} \nonumber\\
 &&\\
 \ope{M}^{\tta_i,\dI\dI\dI}_{j,n}\egd[m^{\tta_i,n}_{k,\ell}]_{\ds^{k\in K(i,j)}_{\ell\in L(i,j)^-}} &\qquad \qquad \qquad& \ope{M}^{\tta_i,\dI\dV}_{j,n}\egd[m^{\tta_i,n}_{k,\ell}]_{\ds^{k\in K(i,j)}_{\ell \in L(i,j)}} \nonumber
 \eeqy
 and then let us define the $\bet_{i,j}\ti \bet_{i,j}$ matrix 
${\bM}^{\tta_i}_{j,n}$   as 
\bbe\la{21021415h2015}\bM^{\tta_i}_{j,n}\egd  \tta_i \left(\ope{M}^{\tta_i,\dI\dV}_{j,n}-\ope{M}^{\tta_i,\dI\dI\dI}_{j,n}\inve{\ope{M}^{\tta_i,\dI}_{j,n}}\ope{M}^{\tta_i,\dI\dI}_{j,n} \right)
\ee
\begin{rem} \label{rem2302142015}
It follows from the fact that the matrix $\bQ$ is invertible,  that $\ope{M}^{\tta_i,\dI}_{j,n}$ is a.s. invertible and so is $\bM^{\tta_i}_{j,n}$. 
\end{rem}


\begin{rem}\la{remherm10162005062015}
In the particular case where $\bA_0$ is Hermitian (which means that $\bQ^{-1} = \bQ^*$ and the $\tta_i$'s are real), then the matrces ${\bM}^{\tta_i}_{j,n}$   are also Hermitian.
\end{rem}

 \indent Now, we can formulate the result on the fluctuations.

\subsubsection{Result}


\begin{theo}\la{belowtheor14122015}
\begin{enumerate}
\item As $N$ goes to infinity, the random vector  $$\displaystyle \left(\Lambda_{i,j,n} \right)_{\substack{1 \leq i \leq q_{\hphantom{1}} \\ 1 \leq j\leq \alp_i \\ 1 \leq n \leq m_i} } $$ defined at \eqre{defLambda30042015} converges    to  the distribution of  a  random vector $$\displaystyle \left(\Lambda^\infty_{i,j,n} \right)_{\substack{1 \leq i \leq q_{\hphantom{1}} \\ 1 \leq j\leq \alp_i \\ 1 \leq n \leq m_i} } $$ with joint distribution defined by the fact that 
for each $1 \leq i \leq q$, $1 \leq j \leq \al_i$ and $1 \leq n \leq m_i$,  $\Lambda_{i,j,n}^\infty$ is the collection of the  ${p_{i,j}}^{\trm{th}}$ roots of the  eigenvalues of some random matrix $\bM^{\tta_i}_{j,n}$. 
\item The distributions of the  random matrices $\bM^{\tta_i}_{j,n}$ are absolutely continuous with respect to the Lebesgue measure  and  the  random vector $\displaystyle \left(\Lambda^\infty_{i,j,n} \right)_{\displaystyle^{1 \leq i \leq q}_{1 \leq j\leq \alp_i}} $   has  no deterministic coordinate. 
\end{enumerate}
\end{theo}

\begin{figure}[ht]
\centering

{
\includegraphics[scale=0.42]{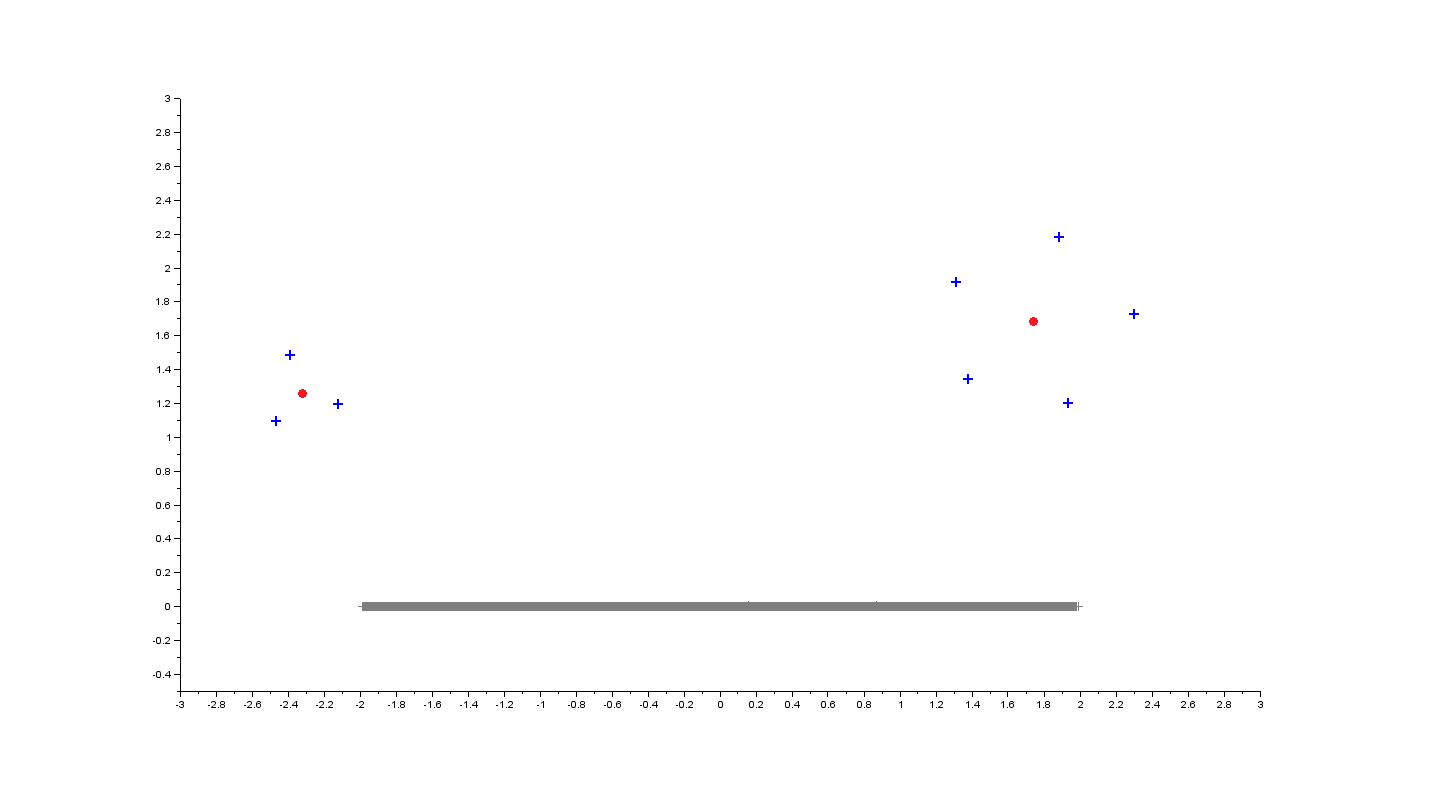}
}
\caption{Spectrum of a Wigner matrix of size $N=5000$ with perturbation matrix $\bA = \diag\lf( \bR_5(1.5+2\ii),\bR_3(-2+1.5\ii),0,\ld,0\ri)$. We   see the blue crosses ``\textcolor{blue}{$+$}'' (outliers) forming respectively a regular pentagon and an equilateral triangle around the red dots ``\tred{$\bullet$}'' (their limit).  We also see a significant difference between the two rates of convergence, $N^{-1/10}$ and $N^{-1/6}$.}\label{wignerfluc1348012015}
\end{figure}

Theorem \ref{belowtheor14122015} is illustrated in Figure \ref{wignerfluc1348012015} with an example. We clearly see appearing regular polygons.

%
%

\section{Applications }

In this section, we give examples of random matrices which satisfy the assumptions of Theorem \ref{theo1123201504058899} and Theorem \ref{belowtheor14122015}. 

\subsection{Wigner matrices} 

Let $\bH_N = \ff{\sqrt N}\bW_N$ be a symmetric/Hermitian Wigner matrix with independent entries up to the symmetry. More precisely, we assume that 
\begin{assum}\la{wigner123510062015} \text{ } \\
\textbf{Real symmetric case :}
\beq
 &&\bullet \ \big( \bW_N\big)_{i,j}, 1 \leq i \leq j \leq N, \ \text{are independent,} \\ 
&& \bullet \ \text{The } (\bW_N)_{i,j}\text{'s} \ \text{ for } i\neq j \text{ (resp. } i=j \text{), are identically distributed,} \qquad\qquad\qquad \qquad \ \\ 
&& \bullet \ \E(\bW_N)_{1,1}=\E(\bW_N)_{1,2} = 0, \ \E(\bW_N)_{1,1}^2 = 2\sigma^2, \ \E(\bW_N)_{1,2}^2 = \sigma^2, \\
&& \bullet \ c_3 := \E\lf| (\bW_N)_{1,1}\ri|^3 < \infty, \ m_5 := \E\lf| (\bW_N)_{1,2}\ri|^5 < \infty. \\
\eeq
\textbf{Hermitian case : } 
\beq &&\bullet \ \big(\re \bW_N\big)_{i,j},\big(\im \bW_N\big)_{i,j}, 1 \leq i < j \leq N, \big(\bW_N\big)_{i,i}, 1 \leq i \leq N, \ \text{are independent}.\\
&& \bullet \ \text{The } (\re\bW_N)_{i,j}\text{'s},(\im\bW_N)_{i,j}\text{'s} \ \text{ for } i\neq j \text{ (resp. } (\bW_N)_{i,i}\text{'s), are identically distributed,}\\ 
&& \bullet \ \E(\bW_N)_{1,1}=\E(\bW_N)_{1,2} = 0, \ \E(\bW_N)_{1,1}^2 = \sigma^2, \ \E(\re\bW_N)_{1,2}^2 = \f{\sigma^2}{2}, \\
&& \bullet \ c_3 := \E\lf| (\bW_N)_{1,1}\ri|^3 < \infty, \ m_5 := \E\lf| (\bW_N)_{1,2}\ri|^5 < \infty. \\
\eeq
\end{assum}
\noindent In this case, we have the following version of Theorem \ref{theo1123201504058899}

\begin{theo}[Convergence of the outliers for Wigner matrices]  \label{wigner183603062015}\text{ } \\
Let $\tta_1,\ld,\tta_j$ be the eigenvalues of $\bA_N$ such that $|\tta_i|>\sigma$. Then, with probability tending to one, for all large enough $N$, there are exactly $j$ eigenvalues $\wt\lam_1,\ld,\wt\lam_j$ of $\wt\bH_N := \ff{\sqrt N}\bW_N + \bA_N$ at a macroscopic distance of $[-2\sigma,2\sigma]$ (outliers). More precisely, for all small enough $\delta>0$, for all large enough $N$, for all $i \in \{1,\ld,j\}$, 
\beq
\wt\lam_{i} & = & \tta_i + \f{\sigma^2}{\tta_i} + \oo1, \ 
\eeq
after a proper labeling. 
\end{theo}

\bpr
We just need to check that Assumptions \ref{assu180203062015} and \ref{assu173602062015} are satisfied. \bgt
\ite[-] As long as the entries of $\bW_N$ have a finite fourth moment, we know (see \cite[Theorem 5.2]{bai-silver-book}) that Assumption \ref{assu180203062015} is satisfied.
\ite[-] Now, we need to show that for any $\delta>0$, as $N$ goes to infinity, 
$$
\sup_{\dist(z,\supp(\mu))>\delta }\norm{\bU_{2r}^* \inve{z\bI - \bH_N} \bU_{2r} - G_{\mu_{sc}}(z) \bI} \ \cvp \ 0.
$$
Since we are dealing with $2r \ti 2r$ sized matrices, it suffices to prove that for any unite vectors $\bu$,$\bv$ of $\C^N$, for any $\delta>0$ and any $\eta>0$, as $N$ goes to infinity, 
$$
\pro\Big(\sup_{\dist(z,\supp(\mu))>\delta }\lf| \bu^* \big(\inve{z\bI - \bH_N}  - G_{\mu_{sc}}(z) \bI \big)\bv\ri| > \eta \Big) \ \tto \ 0.
$$
Moreover, as both $G_{\mu_{sc}}(z)$ and $\norm{\inve{z\bI - \bH_N}}$ goes to $0$ when $|z|$ goes to infinity, we know there is a large enough constant $M$ such that we just need to prove that
$$
\pro\Big(\sup_{\substack{\dist(z,\supp(\mu))>\delta \\ |z| \ \leq \ M}}\lf| \bu^* \big(\inve{z\bI - \bH_N}  - G_{\mu_{sc}}(z) \bI \big)\bv\ri| > \eta \Big) \ \tto \ 0.
$$
Then, for any $\eta'>0$,  the compact set $K = \{z, \ \dist(z,\supp(\mu))>\delta \ \text{ and } \ |z| \leq M \}$ admits a $\eta'$-net, which is a finite set $\{z_1,\ld,z_p\}$ of $K$ such that
$$
\forall z \in K, \ \exists i \in \{1,\ld,p\},  \ \ \ |z-z_i| < \eta',
$$
so that, using the uniform boundedness of the derivative  of $G_{\mu_{sc}}(z)$ and $\bu^*\inve{z - \bH_N}\bv$ on $K$, for a small enough $\eta'$, we just need to prove that 
$$
\pro\Big(\max_{i=1}^p \lf| \bu^* \big(\inve{z_i\bI - \bH_N}  - G_{\mu_{sc}}(z_i) \bI \big)\bv\ri| > \eta/2 \Big) \ \tto \ 0.
$$
Then, we properly decompose each function $ x \ \mapsto \ff{z_i-x}$ as a sum of a smooth compactly supported function and one that vanishes on a neighborhood of $[-2\si,2\si]$ and conclude using \cite[(ii) Theorem 1.6]{PRS11}    
\ent 
Moreover, in the Wigner case, we have
\beq
G_{\mu_{sc}}(z) & = & \f{z - \sqrt{z^2 - 4\sigma^2}}{2\sigma^2},
\eeq
where $\sqrt{z^2 - 4\sigma^2}$ is the branch of the square root with branch cut $[-2\sigma,2\sigma]$  so that for any $z$ outside $[-2\sigma,2\sigma]$, the equation $G_{\mu_{sc}}(z) = \ff \tta$ possesses one solution if and only if $|\tta|> \sigma$ and the unique solution is 
$$
\tta + \f{\sigma^2}{\tta},
$$ 
which means that in the Wigner case, the outliers cannot outnumber the rank of the perturbation, and the \emph{phase transition} condition is simply : $|\tta|>\sigma$. Actually in \cite{BEL} (see Remark 3.2), the authors explain that if $\mu$ is $\boxplus$-infinitely divisible, then the sets $\mathcal{S}_{\tta_i}$'s have at most one element, which means that for Wigner matrices, it is not possible to observe the phenomenon of ``outliers outnumber the rank of $\bA$''.

\begin{rem}
One can find an other proof of Theorem \ref{wigner183603062015} in \cite{ROR1} as a particular case of the Theorem 2.4 (see \cite[Remark 2.5]{ROR1}) due to the fact that a Wigner matrix can be seen as a particular Elliptic matrix. Nevertheless, the authors of \cite{ROR1} don't deal with the matter of the fluctuations. 

\end{rem}

\epr

To study the fluctuations of the outliers in the Wigner case, we must make an additional assumption on the perturbation $\bA_N$.  \\

\begin{assum}\la{ultassu10212052015}

The matrix $\bA_N$ has only a finite number (independent of $N$) of entries which are non-zero.

\end{assum}

\begin{rem}
Assumption \ref{ultassu10212052015} is equivalent to suppose that $\bU_{2r}$ (the $2r$-first columns of $\bU$), possesses only a finite number $K$ (independent of $N$) of non-zero rows. Actually, this assumption is the analogous ``the eigenvectors of $\bA$ don't spread out'' hypothesis corresponding to the ``case a)'' in \cite{CAP}.
\end{rem}

\begin{rem}
If $\bU$ is Haar-distributed and independent from $\bW$, we can avoid making Assumption \ref{ultassu10212052015} (see section \ref{UCIsection20115}). One can also slightly weaken Assumption \ref{ultassu10212052015} by assuming that the $2r$-first rows of $\bU$ correspond to the $N$ first coordinates of a collection of non-random vectors $\bu_1,\ld,\bu_{2r}$ in $\ell^2(\N)$ (see \cite[Theorem 1.7]{PRS11}). 
\end{rem}

\begin{theo}[Fluctuations for Wigner matrices] \text{ } \\
With Assumtions \ref{wigner123510062015} and \ref{ultassu10212052015}, Theorem \ref{belowtheor14122015} holds. Moreover, the distribution of the random vector $$ \lf({m}^{\tta_i}_{k,\ell}\ri)_{\substack{1 \leq i \leq q_{\hphantom{1.}} \quad\quad \ \  \\ (k,\ell)\in J(\tta_i)\ti I(\tta_i)}}, $$ defined by \eqref{2121414h22015}, is
$$
\Big( \be_k^* \bQ^{-1} \bU_{K,2r}^*\Upsilon(\xi_{i})\bU_{K,2r}\bQ   \be_\ell  \Big)_{\substack{1 \leq i \leq q_{\hphantom{1.}} \quad\quad \ \  \\ (k,\ell)\in J(\tta_i)\ti I(\tta_i)}},
$$
where $\xi_i : = \tta_i + \ds\f{\sigma^2}{\tta_i}$ and where $\Upsilon(z)$ is a $K\ti K$ random field defined by
\beqy
\Upsilon(z) & := & (G_{\mu_{sc}}(z))^2 \lf( \bW^{(K)} + \bY(z) \ri)
\eeqy 
where $\bW^{(K)}$ is the $K \ti K$ upper-left corner submatrix of a matrix $\wt\bW_N$ such that $\wt\bW_N \eloi \bW_N$ and $\bY(z)$ is a $K \ti K$ Gaussian random field defined by \cite[(2.7),(2.8),(2.9),(2.10),(2.11),(2.12)]{PRS12} in the real case and \cite[(2.42),(2.43),(2.44),(2.45),(2.46),(2.47)]{PRS12} in the complex case.
\end{theo}

\begin{rem}
This provides   an example of non universal fluctuations, in the sense that the $\lf({m}^{\tta_i}_{k,\ell}\ri)$'s are not necessarily Gaussian. However, when $\bH_N$ is a GOE or GUE matrix, the $\lf({m}^{\tta_i}_{k,\ell}\ri)$'s are centered Gaussian variables such that 
\beqy \la{GOE1005120615}
\E\lf( {m}^{\tta_i}_{k,\ell} \; {m}^{\tta_{i'}}_{k',\ell'}\ri)& =&   \psi_{sc}(\xi_{i},\xi_{i'}) \; \big(\be_k^* \bQ^{-1}(\bQ^{-1})^*\be_{k'} \; \be_\ell^* \bQ^* \bQ \be_{\ell'}\; + \; \delta_{k,\ell'}\delta_{k',\ell}\big), \\ \E\lf( {m}^{\tta_i}_{k,\ell} \; \ovl{{m}^{\tta_{i'}}_{k',\ell'}}\ri) &=&  \psi_{sc}(\xi_{i},\ol{\xi_{i'}}) \; \big(\be_k^* \bQ^{-1}(\bQ^{-1})^*\be_{k'} \; \be_\ell^* \bQ^* \bQ \be_{k'}\; + \; \delta_{k,\ell'}\delta_{k',\ell}\big), \nonumber
\eeqy
for the GOE, and 
\beqy \la{GUE1005120615}
\E\lf( {m}^{\tta_i}_{k,\ell} \; {m}^{\tta_{i'}}_{k',\ell'}\ri)& =&   \psi_{sc}(\xi_{i},\xi_{i'}) \; \delta_{k,\ell'}\delta_{k',\ell}, \\ \E\lf( {m}^{\tta_i}_{k,\ell} \; \ovl{{m}^{\tta_{i'}}_{k',\ell'}}\ri) &=&   \psi_{sc}(\xi_{i},\ol{\xi_{i'}})  \;\be_{ k}^*\bQ^{-1}(\bQ^{-1})^*\be_{ {k'}}\;\be_{ \ell'}^*\bQ^*\bQ\, \be_{ {\ell}}, \nonumber
\eeqy
for the GUE, where
\beq
\psi_{sc}(z,w) & := & G^2_{\mu_{sc}}(z)G^2_{\mu_{sc}}(w)\big(\sigma^2 + \sigma^4  \varphi_{sc}(z,w)\big) , \\
\varphi_{sc}(z,w) &:= &\int \ff{z-x} \ff{w-x} \mu_{sc}(dx).\\
\eeq
 We notice that, if $\bQ^{-1} \neq \bQ^*$, then we might observe correlations between the fluctuations of outliers at a macroscopic distance with each other. 
 This phenomenon has already been observed in \cite{KYIOLW} for non-Gaussian Wigner matrices whereas, here, the phenomenon may still occur for GUE matrices. Actually, \eqref{GOE1005120615} and \eqref{GUE1005120615} can be simplified due to the fact 
 $$
\sigma^2 G_{\mu_{sc}}^2(z) - z G_{\mu_{sc}}(z)+1 \ = \ 0,  
 $$ 
so that  $\varphi_{sc}(z,w)  = \ds-\f{G_{\mu_{sc}}(z) -G_{\mu_{sc}}(w)  }{z-w}$ satisfies
\beqy
\sigma^2 G_{\mu_{sc}}(z)G_{\mu_{sc}}(w) \varphi_{sc}(z,w) & = & \varphi_{sc}(z,w) -  G_{\mu_{sc}}(z)G_{\mu_{sc}}(w). 
\eeqy
Hence, 
\beq
&&G^2_{\mu_{sc}}(\xi_{i})G^2_{\mu_{sc}}(\xi_{i'})\big(\sigma^2 + \sigma^4  \varphi_{sc}(\xi_{i},\xi_{i'})\big) \\& = & 
\sigma^2 G_{\mu_{sc}}(\xi_{i})G_{\mu_{sc}}(\xi_{i'})\Big[ G_{\mu_{sc}}(\xi_{i})G_{\mu_{sc}}(\xi_{i'}) + \sigma^2 G_{\mu_{sc}}(\xi_i)G_{\mu_{sc}}(\xi_{i'}) \varphi_{sc}(\xi_i,\xi_{i'})\Big] \\
& = & \sigma^2 G_{\mu_{sc}}(\xi_{i})G_{\mu_{sc}}(\xi_{i'}) \varphi_{sc}(\xi_i,\xi_{i'}) \\
& = &  \varphi_{sc}(\xi_i,\xi_{i'}) -  G_{\mu_{sc}}(\xi_i)G_{\mu_{sc}}(\xi_{i'}) \ = \ \Phi_{sc}(\xi_i,\xi_{i'}) . 
\eeq
and we fall back on the expression of the variance for the UCI model (see section \ref{UCIsection20115}), which is expected since the GUE belongs to the UCI model.
\end{rem}

\bpr
We show that the assumptions \ref{wigner123510062015} and \ref{ultassu10212052015} imply Assumption \ref{lem15052015}, more precisely $(1)$ and $(2)$. For $(1)$, we simply use \cite[Theorem 2.1/2.5]{PRS12} to show that 
$$
\sqrt N  \bU_{2r}^* \lf( \inve{z-\bH_N} - G_{\mu_{sc}}(z) \bI\ri) \bU_{2r},
$$
converges weakly (as it is done in \cite{PRS11}). The limit distribution is also given by \cite[Theorem 2.1/2.5]{PRS12}. \\
Then for $(2)$, we know by \cite[(i) of Theorem 2.3/2.7]{PRS12} (respectively \cite[(iii) of Proposition 2.1]{PRS12}) that, for all $k \geq 1$, the diagonal entries (respectively the off-diagonal entries) of the matrix $$\sqrt N \big(\lf(z-\bH_N\ri)^{-k-1} - \int (z-x)^{-k-1}\mu_{sc}(dx) \bI\big)$$ converge in distribution so that  
$$
\sqrt N  \bU_{2r}^* \lf( \lf(z-\bH_N\ri)^{-k-1} - \int (z-x)^{-k-1}\mu_{sc}(dx) \bI\ri) \bU_{2r}
$$
is tight.
\epr

\subsection{Hermitian matrices whose distribution is invariant by unitary conjugation} \label{UCIsection20115}
Let $\bH_N$ be an Hermitian matrix such that for any unitary $N \ti N$ matrix $\bU_N$, we have
\beqy
\bU_N \bH_N \bU_N^* & \eloi & \bH_N.
\eeqy
$\bH_N$ can be written $\bH_N = \bU_N \bD_N \bU_N^*$ where $\bD_N$ is diagonal, $\bU_N$ is Haar-distributed and $\bU_N$ and $\bD_N$ are independent. We also assume that $\bH_N$ satisfies \eqref{assu1055m2015k} and Assumption \ref{assu180203062015}. We shall call such matrices \emph{UCI matrices} (for Unitary Conjugation Invariance). In this case, as we can  we can write
$$
\wt\bH_N \ = \ \bH_N + \bA_N \ = \ \bU_N \lf( \bD_N + \bU^*_N \bA_N \bU_N \ri) \bU_N^*, 
$$
so that, without any loss of generality, we can simply assume that $\bH_N$ is a diagonal matrix and $\bA_N$ is a matrix of the form 
\beq
\bA_N & = & \bU_{2r} \bA_0 \bU_{2r}^*
\eeq
where $\bU_{2r}$ is the $2r$-first columns of an Haar-distributed matrix independent from $\bH_N$.

\begin{theo}[Convergence of the outliers for UCI matrices] \la{theo121211123201504058899} \text{ } \\
If $\bH_N$ is an UCI matrix, then Theorem \ref{theo1123201504058899} holds. 
\end{theo}

\begin{rem}
Unlike the Wigner case, Theorem \ref{theo1123201504058899} does not need to be reformulated. In this case, we do observe the phenomenon of the outliers outnumbering the rank of $\bA_N$. 
\end{rem}

\bpr
We just need to check that Assumption \ref{assu173602062015} is satisfied. To do so, one can apply a slightly modified version of \cite[Lemma 2.2]{BEN2}, where we replace all the ``$\dist(z,[a,b])>\delta$'' by ``$\dist(z,\supp(\mu))>\delta$'', which does not change the ideas of the proof. 
\epr

For the fluctuations, we need to assume that for all $i\in \{1,\ld,q\}$ and all $n\in \{1,\ld,m_i\}$, as $N$ goes to infinity, 
\beqy
\sqrt{N}\lf(\ff N \tr \inve{\xi_{i,n}-\bH_N} - \ff{\tta_i}\big)\ri) \tto 0.
\eeqy

\begin{rem}
Actually, in \cite{BEN2}, the authors make the same assumption (\cite[Hypothesis 3.1]{BEN2}).
\end{rem}

\begin{theo}[Fluctuations for UCI matrices] \text{ } \\
If $\bH_N$ is an UCI matrix, then if satisfies Theorem \ref{belowtheor14122015}. More precisely, the  $$ \lf({m}^{\tta_i,n}_{k,\ell}\ri)_{\substack{1 \leq i \leq q_{\hphantom{1.}} \quad\quad \ \ \\ 1 \leq n \leq m_i \quad\quad \ \ \\ (k,\ell)\in J(\tta_i)\ti I(\tta_i)}}, $$ defined by \eqref{2121414h22015} are centered Gaussian variables such that 
\beq
\E\lf( {m}^{\tta_i,n}_{k,\ell} \; {m}^{\tta_{i'},n'}_{k',\ell'}\ri)& =&   \Phi(\xi_{i,n},\xi_{i',n'}) \; \delta_{k,\ell'}\delta_{k',\ell}, \\ \E\lf( {m}^{\tta_i,n'}_{k,\ell} \; \ovl{{m}^{\tta_{i'},n'}_{k',\ell'}}\ri) &=&  \Phi(\xi_{i,n},\ol{\xi_{i',n'}}) \;\be_{ k}^*\bQ^{-1}(\bQ^{-1})^*\be_{ {k'}}\;\be_{ \ell'}^*\bQ^*\bQ\, \be_{ {\ell}},
\eeq
where 
\beq
\Phi(z,w) & := & \int \ff{z-x}\ff{w-x}\mu(dx) - \int \ff{z-x}\mu(dx)\int \ff{w-x}\mu(dx)\\ & = & \begin{cases} - \f{G_\mu(z)-G_\mu(w)}{z-w} - G_\mu(z)G_\mu(w)& :  \text{ if } z \neq w, \\
                                                                     - G'_\mu(z) - (G_\mu(z))^2& : \text{ otherwise.} \\ \end{cases} \\
\eeq
\end{theo}

\begin{rem}
Remind that we supposed that $\mu$ is not a single Dirac measure, so that $\Phi$ is not equal to zero. 
\end{rem}

\begin{figure}[ht]
\centering
\subfigure[Uncorrelated case : $G_\mu(\xi_1)\neq G_\mu(\xi_2)$, which means that $\xi_1$ and $\xi_2$ do not belong to the same set $\mathcal{S}_\tta$.]{\includegraphics[scale=0.3]{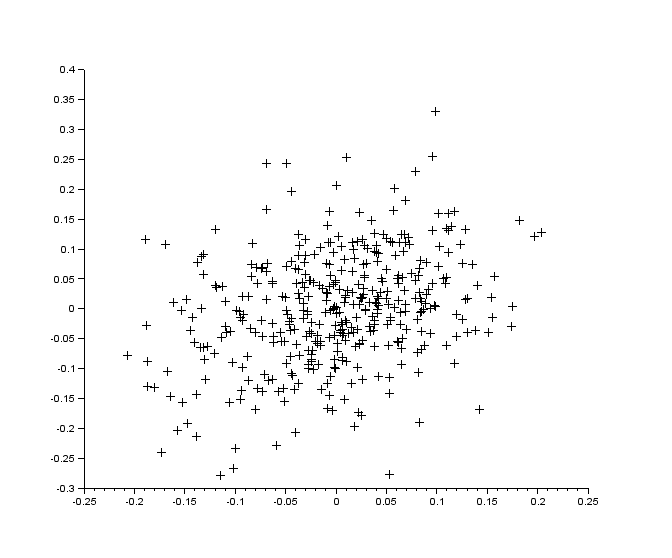}\label{contrexemple14352220150011445bis}}\qquad
\subfigure[Correlated case : $G_\mu(\xi_1)= G_\mu(\xi_2)$, which means that $\xi_1$ and $\xi_2$ belong to the same set $\mathcal{S}_\tta$.]{\includegraphics[scale=0.3]{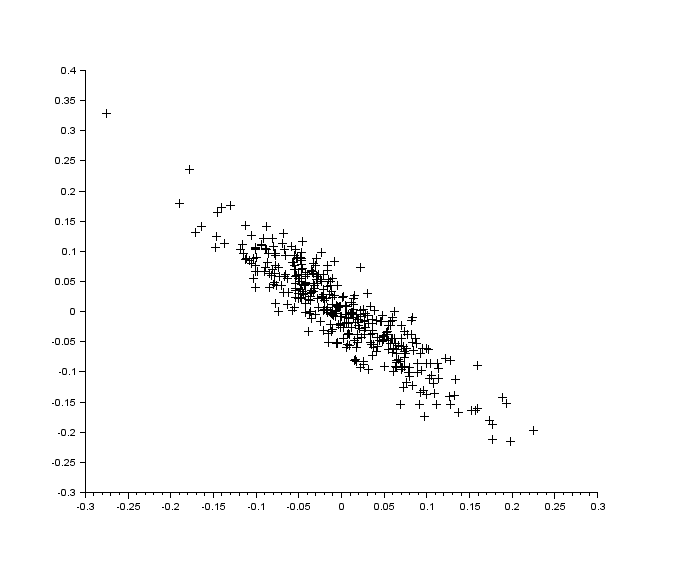} \label{contrexemple14352220150011445}}
\caption{Correlation between the fluctuations of two outliers $\xi_1\neq\xi_2$ for a sample of 400 matrices of size 500, in the case where both $\bH$ and $\bA$ are Hermitian and where the support of $\mu$ is disconnected. Here $\mu(dx)$ is taken equal to $\ff 2(\delta_{-1}(dx)+\one_{[1,2]}(dx))$.}\label{refFigure17154567}
\end{figure}

\begin{rem}
If $\bA_N$ is Hermitian, the size of all the Jordan blocks are equal to $1$ and the fluctuations are real random variables (see Remark \ref{remherm10162005062015}). We find back that, in the Hermitian case, fluctuations between outliers at a macroscopic distance are independent (see \cite{BEN2}) except if the two outliers come from the same eigenvalue of $\bA$ (i.e. they both belong to the same set $\mathcal{S}_\tta$).   In this case, the fluctuations of outliers belonging to the same set $\mathcal{S}_\tta$ are all correlated.  This phenomenon is illustrated by Figures \ref{contrexemple14352220150011445bis} and \ref{contrexemple14352220150011445}. 
\end{rem}

\bpr We just need to check that $\bH_N$ satisfies $(1'),(2')$ of Assumption \ref{lem15052015} (since $(0')$ is assumed below). Actually, for any $k \geq 1$ and any $i \in \{1,\ld,q\}$, the diagonal matrix
$$
\big(\xi_{i,n}-\bH_N\big)^{-(k+1)} -  \ff N \tr\big(\xi_{i,n}-\bH_N\big)^{-(k+1)},
$$  
fulfill the assumptions of Theorem \ref{lem095815052015}, so that $(2')$ is true. Then, $(1')$ is true thanks to Theorem \ref{theo095815052015}. This theorem also gives us the covariance.

\epr

\section{Proofs}

\subsection{Convergence of the outliers : proof of Theorem \ref{theo1123201504058899}}

In \cite{BEL}, the authors give an interpretation of why the limit is necessarily a solution of $G_{\mu}(z) = \ff \tta$ with the subordinate functions of the free additive convolution of measures in the particular case where one of the measure is $\delta_0$ (see \cite[Example 4.1]{BEL}). Actually, our definition of the sets $\mathcal{S}_{\tta_i}$'s corresponds to the one of the set $O_\tta$ in \cite[Definition 4.1]{BEL}. A quick (but inaccurate) way to see why the limit is  $G_{\mu}^{-1}(\ff \tta)$ and to understand the approach of the proof, is to write
$$
\det \lf(z - (\bH_N+\bA_N) \ri) \ = \ \det \lf( z - \bH_N \ri)\det \lf( \bI - \inve{z - \bH_N}\bA_N \ri),
$$
then if $\inve{z - \bH_N} \sim G_\mu(z)\bI$, we can write
$$
\det \lf(z - (\bH_N+\bA_N) \ri) \sim \det \lf( z - \bH_N \ri) G_\mu(z) \det \lf( \ff{G_\mu(z)}\bI - \bA_N \ri)
$$
so that if $z$ is an outlier of $\bH_N+\bA_N$, $\ff{G_\mu(z)}$ must be an eigenvalue of $\bA_N$.\\

To do it properly, we introduce the following function\footnote{we used a classical trick of finite rank perturbation models which $\det(\bI_m+\bA\bB)= \det(\bI_n + \bB\bA)$ for any $m \ti n$ matrix $\bA$ and $n \ti m$ matrix $\bB$},
\beqy
f(z) & = & \det\left(\bI - \bU_{2r}^*\inve{z\bI - \bH}\bU_{2r}\bA_0 \right) \ = \ \f{\det \lf(z - (\bH_N+\bA_N) \ri) }{\det \lf( z - \bH_N \ri)},
\eeqy
we know that the zeros of $f$ are eigenvalues of $\wt\bH_N$ which are not eigenvalues of $\bH_N$. Then, we introduce the function
\beqy
f_0(z) &:= & \det \left( \bI - G_\mu(z) \bA_0\right),
\eeqy
and the proof of Theorem \ref{theo1123201504058899} relies on the two following lemmas.
\begin{lemme} \la{lem1452060502015}
As $N$ goes to infinity, we have
$$
\sup_{\dist(z,\supp(\mu))>\delta } \lf|f(z)-f_0(z)\ri| \ \cvp \ 0.
$$
\end{lemme}
\begin{lemme} \label{lem1131020505201055}
Let $K$ be a compact set and let $\ep>0$ such that \bgt
\ite $\ds\dist(K,\bigcup_{i=1}^j \mathcal{S}_{\tta_i})\geq \ep$, 
\ite $\ds\dist(K,\supp(\mu)) \geq \ep$.
\ent
Then, with a probability tending to one, 
$$
\inf_{z \in K}\lf|\det \lf(\bI - \inve{z - \bH_N} \bA_N\ri)\ri| \ > \ 0.
$$

\end{lemme}

If these lemmas are true, the end of the proof goes as follow. We know that, with a probability tending to one, there is $\ep>0$, such that \bgt
\ite there is a constant $M>0$ such that $\bH_N+\bA_N$ has no eigenvalues in the area $\{z, \ |z|>M\}$,
\ite  $\spec(\bH_N)\subset\{z, \ \dist(z,\supp(\mu))<\ep\}$,
\ent
We set 
$$
\mathcal S \ := \ \bigcup_{i=1}^j \mathcal{S}_{\tta_i}
$$
and we define 
\beqy
\mathcal{S}^\ep & := & \bigcup_{i=1}^j \bigcup_{\xi \in \mathcal{S}_{\tta_i}} \lf\{z, \ |z-\xi|<\ep\ri\} 
\eeqy
with the convention that ${\mathcal S}^{\ep} = \emptyset$ if $\mathcal S=\emptyset$. Up to a smaller choice of $\ep$, we can suppose that none of the disk centered in the element of the $\mathcal{S}_{\tta_i}$'s and of radius $\ep$ intersects each other nor intersect $\{z, \ \dist(z,\supp(\mu))<\ep\}$. Then, using Lemma \ref{lem1131020505201055}, with
$$
K \ := \ \lf\{ z, \ |z| \leq M \ri\} \backslash \lf( {\mathcal S}^{\ep} \cup \{z, \ \dist(z,\supp(\mu))<\ep\} \ri),
$$
we deduce all the eigenvalues of $\wt \bH_N$ are contained in ${\mathcal S}^{\ep } \cup \{z, \ \dist(z,\supp(\mu))<\ep\}$. Indeed, if $z$ is an eigenvalue of $\wt\bH_N$ such that $\dist(z,\supp(\mu))>\ep$, $z$ must be a zero of $f$. \\ Moreover, for each $i\in\{1,\ld,j\}$ and each $\xi \in \mathcal{S}_{\tta_i}$, we know that from Lemma \ref{lem1452060502015}
\beq
\sup_{z, |z-\xi|=\ep} |f(z)-f_0(z)| \ \tto \ 0, & \quad \text{ and } \quad & \inf_{z, |z-\xi|=\ep} |f_0(z)| \ > \ 0.
\eeq
 we deduce by Rouch\'e Theorem (see \cite[p. 131]{Beardon}) that $f$ and $f_0$, for all large enough $N$, have the same number of zeros inside the domain $\{z, \ |z-\xi|<\ep\}$, for each $\xi$ in the $\mathcal{S}_{\tta_i}$'s. \\
 Now, we just need to prove the two previous lemmas.

\bpr[of Lemma \ref{lem1452060502015}]
We know that, for some positive constant $C$,
\beq
\sup_{\dist(z,\supp(\mu))>\delta }|f(z) - f_0(z)| & \leq & C \sup_{\dist(z,\supp(\mu))>\delta }\norm{\bU_{2r}^*\lf( \inve{z - \bH_N}-G_\mu(z)\bI\ri)\bU_{2r}}
\eeq
and we conclude with Assumption \ref{assu173602062015}.
\epr

\bpr[of Lemma \ref{lem1131020505201055}]
We write, thanks to Assumption \ref{assu173602062015},
\beq
\det \lf( \bI - \inve{z - \bH_N} \bA_N\ri) & = & \det \lf(\bI - \bU_{2r}^* \inve{z-\bH_N}\bU_{2r}\bA_0 \ri) \\
& = & \det \lf( \bI - G_{\mu}(z) \bA_0 +\oo1\ri) \\
& = & \prod_{i=1}^k \big(1 - G_{\mu}(z) \tta_i \big)+\oo1 .
\eeq
Then, since $z \in K$, it easy to show that for each $i$, $|1 - G_{\mu}(z)\tta_i|>0$.
\epr
%

\subsection{Fluctuations}
The proof of Theorem \ref{belowtheor14122015} is the same than \cite[Theorem 2.10]{FloJean} and all we need to do here is to prove this analogous version of \cite[Lemma 5.1]{FloJean}. 
\begin{lemme} \label{lem112107010101017}
\indent For all $j \in \{1,\ldots,\alp_i\}$ and all $n \in \{1,\ld,m_i\}$, let $F^{\tta_i}_{j,n}(z)$ be the rationnal function defined by
 \bbe\la{162142bis}
F^{\tta_i}_{j,n}(z) \ : = \ f\left( \xi_{i,n} + \frac{z}{N^{1/(2p_{i,j})}} \right).
\ee
\indent Then, there exists a collection  of  positive constants $(\gamma_{i,j})_{\ds^{1\le i\le q}_{1\le j\le \al_i}}$ and a collection of non vanishing random variables $(C_{i,j,n})_{\substack{1\le i\le q \\ 1\le j\le \al_i \\ 1 \le n \le m_i}}$ independent of $z$,   such that we have the convergence in distribution (for the topology of the uniform convergence over any compact set)
$$
\lf(  N^{\gamma_{i,j}} F^{\tta_i}_{j,n}(\cdot) \ri)_{\substack{1\le i\le q \\ 1\le j\le \al_i \\ 1 \le n \le m_i}}\ \underset{N \to \infty}{\tto} \  \left( z\in \C \ \mapsto \ z^{\pi_{i,j}}\cdot  C_{i,j,n} \cdot \det\left( z^{p_{i,j}}-\bM^{\tta_i}_{j,n}\right) \right)_{\substack{1\le i\le q \\ 1\le j\le \al_i \\ 1 \le n \le m_i}}
$$
    where $\bM^{\tta_i}_{j,n}$ is the random matrix introduced at \eqref{2121414h32015} and  $\pi_{i,j} \egd \sum_{l > j} \bet_{{i,l}}p_{i,l}$. 
\end{lemme}
Once this lemma proven, the Theorem \ref{belowtheor14122015} follows (see section 5.1 of \cite{FloJean} for more details). To prove Lemma \ref{lem112107010101017}, we shall proceed as it is done in \cite{FloJean} to prove Lemma 5.1. First, we write, for a fixed $\tta_i(=\tta)$, a fixed $n\in\{1,\ld,m_i\}$ and a fixed $j\in\{1,\ld,\al_i\}$ (which shall be implicit) and fixed $p_{i,j}$($=p$), recall that $\bA_0 = \bQ \bJ \bQ^{-1}$, 
\beq
F^\tta_{j,n}(z) & =& \det \lf( \bI - \inve{\xi_n + \f{z}{N^{1/(2p)}}-\bH_N} \bU_{2r}\bQ\bJ\bQ^{-1}\bU_{2r}^*\ri)\\
 & = & \det \left( \bI - G_\mu\Big( \xi_n+\f{z}{N^{1/(2p)}}\Big)\bJ - \ff{\sqrt{N}}{ \bZ_N\Big( \xi_n + \f{z}{N^{1/(2p)}} \Big)}\right) \\
& = & \det \left( \bI -\f{\bJ}{\tta}-G'_\mu(\xi_n)\f{z}{N^{1/(2p)}}\big(1+\oo{1}\big)\bJ  -\ff{\sqrt{N}}{ \bZ_N\Big(\xi_n +\f{z}{ N^{1/(2p)}} \Big)}\right) \\
\eeq
where
$$
Z_N(z) \ := \ \sqrt{N} \bQ^{-1}\bU_{2r}^*\lf( \inve{z-\bH_N}-G_\mu(z)\bI\ri)\bU_{2r}\bQ \bJ.
$$ 
Remind that by definition, $G_\mu(\xi_n) = \tta^{-1}$. From here, the reasoning to end the proof is the exact same than the one from \cite[Lemma 5.1]{FloJean}. Nevertheless, we still have to prove that, for all $\tta$ and for all $n$, for all compact set $K$ and for all $z \in K$, 
\beqy \label{equ102613052015}
\bZ_N\big(\xi_n+\f{z}{N^{1/(2p)}}\big) & = & \bZ_N(\xi_n) + \oo1, \qquad \text{ and   $\ \ \bZ_N(\xi_n)$ converges weakly. }
\eeqy
To do so, we write (thanks to \ref{lemm150020151205}),
\beq
 \bZ_N\lf( \xi_n + \f{z}{N^{1/(2p)}}\ri)
& = & \sqrt{N}\bQ^{-1}\bU_{2r}^*\lf( \inve{\xi_n-\bH_N} - \ff{\tta}\ri)\bU_{2r}\bQ\bJ \\ && \ \ + \ \sum_{k=1}^p  \lf(\f{-z}{N^{1/(2p)}}\ri)^k \sqrt{N}\bQ^{-1}\bU_{2r}^* \lf(\big(\xi_n-\bH_N\big)^{-(k+1)} - \int \f{\mu(dx)}{(\xi_n-x)^{k+1}} \ri)\bU_{2r}\bQ\bJ   \\
&  & \ \  + \ \ff{N^{1/(2p)}} \bQ\bU_{2r}^*\Big({\xi_n}-\bH_N \Big)^{-(k+1)}\inve{\xi_n+\f{z}{N^{1/(2p)}} -\bH_N}\bU_{2r}\bQ\bJ + \oo1. 
\eeq
The last term is a $\oo1$ since $\dist(\xi_n,\spec(\bH_N))> \ep$ and one can conclude if $(1),(2)$ are satisfied in Assumption \ref{lem15052015}. Otherwise, if it's $(0'),(1'),(2')$, we write
\beq
 \bZ_N\lf( \xi_n + \f{z}{N^{1/(2p)}}\ri) & = & \sqrt N \bQ^{-1}\bU_{2r}^* \lf( \ff N \tr \inve{\xi_n-\bH_N} - \ff {\tta} \ri) \bU_{2r}\bQ\bJ \\
 &&\  + \ \sqrt{N}\bQ^{-1}\bU_{2r}^*\lf( \inve{\xi_n-\bH_N} - \ff N \tr \inve{\xi_n-\bH_N} \ri) \bU_{2r}\bQ\bJ  \\
 &&\  + \ \sum_{k=1}^p  \lf( \f{-z}{N^{1/(2p)}} \ri)^k \sqrt{N}\bQ^{-1}\bU_{2r}^* \lf(\big(\xi_n-\bH_N\big)^{-(k+1)} - \ff N \tr \inve{\xi_n-\bH_N} \ri)\bU_{2r}\bQ\bJ   \\
&  &  \  + \ \ff{N^{1/(2p)}} \bQ\bU_{2r}^*\Big({\xi_n}-\bH_N \Big)^{-(p+1)}\inve{\xi_n+\f{z}{N^{1/(2p)}} -\bH_N}\bU_{2r}\bQ\bJ + \oo1. 
\eeq

\appendix

\section{}
\subsection{Linear algebra lemmas}
\begin{lemme}\la{lemm150020151205}
Let $\bA$ be a matrix and $\lam \in \C$ be such that both $\bA$ and $\bA+\lam\bI$ are non singular. Then, for all $p \geq 1$,
\beq
\inve{\bA + \lam\bI} & = & \sum_{k=1}^p (-\lam)^{k-1}\bA^{-k} + (-\lam)^{p} \bA^{-p}\inve{\bA+\lam\bI} 
\eeq
\end{lemme}

\begin{lemme}[Schur's complement  \cite{HORNJOHNSON} ] \la{schurcom15222015} For any $\bA,\bB,\bC,\bD$, one has, when it makes sense
$$
\bpm \bA & \bB \\ \bC & \bD \\ \epm^{-1} \ = \ \bpm \inve{\bA - \bB \bD^{-1} \bC} & -\inve{\bA - \bB \bD^{-1} \bC}\bB\bD^{-1} \\
- \bD^{-1}\bC\inve{\bA - \bB \bD^{-1} \bC} & \bD^{-1} + \bD^{-1}\bC\inve{\bA - \bB \bD^{-1} \bC}\bB\bD^{-1} \epm
$$
\end{lemme}

\subsection{Fluctuations of the entries of UCI random matrices}

We give here some results on the fluctuations of the entries of UCI matrices, which means, matrices of the form $\bH:=\bU\bD\bU^*$ where $\bU$ is Haar-distributed and $\bD$ is a complex diagonal matrix.

\begin{theo}[Fluctuations of the entries of UCI random matrices]\la{lem095815052015}
Let $\bT$ be an $N \ti N$ diagonal matrix such that 
\beqy \la{09540205062015}
\qquad\tr \bT \ = \ 0, \quad \quad \ff N \tr \bT \bT^* \ \to \ \sigma^2, \quad \quad \ff N \tr \bT^2 \ \to \ \tau^2, \quad \quad \forall k \geq 1, \ \ff N \tr (\bT \bT^*)^k \ = \ \OO1.
\eeqy
Let $\bu_{t_1},\ld,\bu_{t_p}$ be   $p$ distinct columns of a Haar-distributed unitary matrix. Then
$$
\lf(\sqrt N \scalv{\bu_{t_i}}{\bT\bu_{t_j}}\ri)^{p}_{i,j=1},
$$
converges in distribution to a centered complex  Gaussian vector $\big(\mathcal{G}_{i,j}\big)_{i,j=1}^p$ with covariance
\beq
\Ec{\mathcal{G}_{i,j}\mathcal{G}_{k,\ell}} \ = \   \delta_{i,\ell}\delta_{j,k} \tau^2         &;&
\Ec{\mathcal{G}_{i,j}\ol{\mathcal{G}}_{k,\ell}} \ = \ \delta_{i,k}\delta_{j,\ell}\sigma^2 
\eeq
\end{theo}
\begin{rem}
If $\bH := \bU \bD \bU^*$ satisfies \eqref{assu1055m2015k} and Assumption \ref{assu180203062015}, then $\bT:=\bD-\ff N \tr \bD$ satisfies \eqref{09540205062015}.
\end{rem}

Here comes a version of Theorem \ref{lem095815052015}, with several matrices diagonal $\bT$. Due to the complex values of the diagonal matrices, the following theorem is not a simple consequence of Theorem \ref{lem095815052015} and Cram\'er--Wold theorem.


\begin{theo}\la{theo095815052015}
Let $\bT_1,\ld,\bT_q$ be $N \ti N$ diagonal matrices such that for all $m, n \in \{1,\ld,q\}$ 
\beq
\tr \bT_m \ = \ 0, \quad \quad \ff N \tr \bT_m \bT^*_n \ \to \ \sigma^2_{m,n}, \quad \quad \ff N \tr \bT_m\bT_n \ \to \ \tau^2_{m,n}, \quad \quad \forall k \geq 1, \ \ff N \tr (\bT_m \bT^*_m)^k \ = \ \OO1.
\eeq
Let $\bu_{t_1},\ld,\bu_{t_p}$ be   $p$  distinct columns of an Haar-distributed matrix. Then
$$
\lf(\sqrt N \scalv{\bu_{t_i}}{\bT_m\bu_{t_j}}\ri)_{\substack{1 \leq i \leq p \\ 1 \leq j \leq p \\ 1 \leq m \leq q } },
$$
converges in distribution to a centered complex  Gaussian vector $\big(\mathcal{G}_{i,j,m}\big)_{\substack{1 \leq i \leq p \\ 1 \leq j \leq p \\ 1 \leq m \leq q } }$ with covariance
\beq
\Ec{\mathcal{G}_{i,j,m}\mathcal{G}_{k,\ell,n}} \ = \   \delta_{i,\ell}\delta_{j,k} \tau^2_{m,n}         &;&
\Ec{\mathcal{G}_{i,j,m}\ol{\mathcal{G}}_{k,\ell,n}} \ = \ \delta_{i,k}\delta_{j,\ell}\sigma^2_{m,n} 
\eeq
\end{theo}

\noindent\textbf{Proof of Theorem \ref{lem095815052015}.} Without any loss of generality, due to the invariance by conjugation by a matrix of permutation, we can suppose that $t_1=1,t_2=2,\ld,t_p=p$. Then, we just need to show that 
$$
X \ = \ \sqrt N \tr\big( \bU^* \bT \bU \bA\big)
$$
where $\bA$ is a $N\ti N$ deterministic matrix of the form
$$
\bA \ = \ \bpm \bA_p & 0 \\ 0 & 0 \epm, \ \quad \bA_p = \big( a_{i,j}\big)_{i,j=1}^p \in \M_p(\C),
$$
is a asymptotically Gaussian. Before starting, we remind some definition. Let $(\bM_1,\ld,\bM_q)$ be $q$ matrices. 
For any permutation $\si \in S_q$, with cycle decomposition 
$$
\si \ = \ (i_{1,1} \cdots i_{1,k_1})(i_{2,1}\cdots i_{2,k_2}) \cdots (i_{r,1}\cdots  i_{r,k_r})
$$
 we denote by
\beqy \la{hyp0949201520151}
\tr_{\sigma} \big( \bM_t \big)_{t=1}^q & : = & \prod_{j=1}^r \tr\big( \bM_{t_{i_{j,1}}} \cdots \bM_{t_{i_{j,r_j}}}\big).
\eeqy
For example, if $\sigma = (13)(256) \in S_6$, then
$$
\tr_\si \big(\bM_t \big)_{t=1}^6 \ = \ \tr (\bM_1 \bM_3) \tr (\bM_2\bM_5\bM_6) \tr (\bM_4).
$$
Let $M(2n)$ be the set of all \emph{perfect matching} on $\{1,\ld,2n\}$ which is a subset of $S_{2n}$ of the permutation which are the product of $n$ transpositions with disjoint support. For example
$$
M(4) \ = \ \lf\{(12)(34),(13)(24),(14)(23)\ri\}.
$$
Then, if the following lemma is true, one can conclude the proof.
\begin{lemme}\la{lem151727052015}
Let $\bT_1 ,\ld, \bT_q$ be $q$ diagonal matrix such for all $i,j \in \{1,\ld,q\}$,
\beqy
\tr \bT_i \ = \ 0 \quad  ; \quad \ff N \tr \bT_i \bT_j \ \tto \ \tau_{i,j}  \quad ; \quad \forall k\geq 1, \ \ff N \tr \big( \bT_i \bT_i^*\big)^k \ = \ \OO1 .
\eeqy
Let $\bA_1,\ld,\bA_q$ be $q$ matrices of the form
$$
\bA_i \ = \ \bpm \bA_{0,i} & 0 \\ 0 & 0 \\ \epm 
$$
where the $\bA_{0,i}$'s are $K \ti K$ matrices independent from $N$ where $K$ is a fixed integer. Let $\bU$ be a Haar-distributed matrix. Then, as $N$ goes to infinity,
\beq
\Ec{\prod_{t=1}^q \sqrt N \tr\big( \bU^* \bT_t \bU \bA_t\big)} & \tto & \begin{cases}\ds \sum_{\sigma \in M(q)}  \tr_\sigma\big(\bA_i  \big)_{i=1}^q \prod_{t=1}^{q/2}\tau_{\sigma(2t-1),\sigma(2t)} & \text{if } q \text{ is even.}\\
\quad \qquad 0 & \text{if } $q$ \text{ is odd.}\end{cases}
\eeq
\end{lemme}

Indeed, once we suppose Lemma \ref{lem151727052015} satisfied, we need to compute for all $p,q$
$$
\Ec{\big[\sqrt{N}\tr \bU^* \bT \bU \bA\big]^p \big[\ol{\sqrt N \tr\bU^* \bT \bU \bA}\big]^q} \ = \ \Ec{\big[\sqrt{N}\tr \bU^* \bT \bU \bA\big]^p \big[{\sqrt N \tr\bU^* \bT^* \bU \bA^*}\big]^q} ,
$$
in order to apply Lemma \ref{lemgauss15182015}. According to Lemma \ref{lem151727052015}, for $\bT_t \equiv \bT$ and $\bA_t \equiv \bA$, we have
\beq
\E\big[\sqrt{N}\tr \bU^* \bT \bU \bA\big]^q & \tto & \begin{cases} \card M(q)  \tr(\bA^2)^{q/2}  \tau^{q/2} & \text{if } $q$ \text{ is even,}\\
0& \text{if } $q$ \text{ is odd.}
\end{cases} 
\eeq
(remind that $\card M(q) = (q-1)(q-3)\cdots 3$) which means that the limit distribution of $X$ already satisfies \eqref{cond11444201505} and \eqref{cond21444201505}. Let $p \geq 1$ and $q \geq 2$ be two fixed integers such that $p+q$ is even, then, using notations from \eqref{hyp0949201520151}, we know thanks to Lemma \ref{lem151727052015} that

\beqy\la{eq155905052720155}
\qquad\Ec{\big[\sqrt{N}\tr \bU^* \bT \bU \bA\big]^p \big[{\sqrt N \tr\bU^* \bT^* \bU \bA^*}\big]^q} & = & \ff{N^{\f{p+q}{2}}} \sum_{\sigma \in M(p+q)}  \tr_\sigma\big(\bA_t  \big)_{t=1}^{p+q}\tr_\sigma\big(\bT_t\big)_{t=1}^{p+q} + \oo1
\eeqy
where $$(\bT_1,\ld,\bT_{p+q}) \ = \ (\underbrace{\bT,\ld,\bT}_{p},\underbrace{\bT^*,\ld,\bT^*}_{q}) \ \ \text{ and } \ \ (\bA_1,\ld,\bA_{p+q}) \ = \ (\underbrace{\bA,\ld,\bA}_{p},\underbrace{\bA^*,\ld,\bA^*}_{q}).
$$
We rewrite the right side of \eqref{eq155905052720155} summing according to the value of $\sigma(1)$.
\beq
\sum_{\sigma \in M(p+q)}  \tr_\sigma\big(\bA_t  \big)_{t=1}^{p+q}\tr_\sigma\big(\bT_t\big)_{t=1}^{p+q}  & = & \sum_{a=2}^{p+q}\sum_{\substack{\sigma \in M(p+q) \\ \sigma(1)=a}}  \tr_\sigma\big(\bA_t  \big)_{t=1}^{p+q}\tr_\sigma\big(\bT_t\big)_{t=1}^{p+q} \\
& = &  \sum_{a=2}^{p+q}\tr(\bA_1 \bA_a) \tr(\bT_1 \bT_a) \sum_{\substack{\sigma \in M(p+q) \\ \sigma(1)=a}}  \tr_{\sigma\circ(1a)}\big(\bA_t  \big)_{t=1}^{p+q}\tr_{\sigma\circ(1a)}\big(\bT_t\big)_{t=1}^{p+q} \\
& = & \sum_{a=2}^p \tr\bA^2 \tr \bT^2  \sum_{\sigma \in M(p+q-2) }  \tr_{\sigma}\big(\widehat\bA_t  \big)_{t=1}^{p+q-2}\tr_{\sigma}\big(\widehat\bT_t\big)_{t=1}^{p+q-2} \\
&  & \qquad +  \sum_{a=p+1}^{p+q} \tr\bA\bA^* \tr \bT\bT^*  \sum_{\sigma \in M(p+q-2) }  \tr_{\sigma}\big(\wt\bA_t  \big)_{t=1}^{p+q-2}\tr_{\sigma}\big(\wt\bT_t\big)_{t=1}^{p+q-2},
\eeq
where $$(\widehat\bA_1,\ld,\widehat\bA_{p+q-2}) \ = \ (\underbrace{\bA,\ld,\bA}_{p-2},\underbrace{\bA^*,\ld,\bA^*}_{q}) \ \ \text{ and } \ \ (\wt\bA_1,\ld,\wt\bA_{p+q-2}) \ = \ (\underbrace{\bA,\ld,\bA}_{p-1},\underbrace{\bA^*,\ld,\bA^*}_{q-1}).
$$
At last, one easily deduces that
\beq
\hspace{-2.5mm}\Ec{\big[\sqrt{N}\tr \bU^* \bT \bU \bA\big]^p \big[{\sqrt N \tr\bU^* \bT^* \bU \bA^*}\big]^q}   =  \ff N \tr \bT^2 \tr \bA^2 (p-1)\Ec{\big[\sqrt{N}\tr \bU^* \bT \bU \bA\big]^{p-2} \big[{\sqrt N \tr\bU^* \bT^* \bU \bA^*}\big]^q} \\
 \ + \ \ff N \tr \bT\bT^* \tr \bA\bA^* q\Ec{\big[\sqrt{N}\tr \bU^* \bT \bU \bA\big]^{p-1} \big[{\sqrt N \tr\bU^* \bT^* \bU \bA^*}\big]^{q-1}}+\oo1
\eeq
and so $\sqrt N \tr (\bU^* \bT \bA \bU)$ satisfies \eqref{cond31444201505} which means according to Lemma \ref{lemgauss15182015} that its limit distribution is Gaussian. \\
At last, to compute to covariance of the $\big(\mathcal{G}_{i,j}\big)$'s, one can simply use \cite[Lemma A.6]{FloJean2}.
\hfill$\square$ \\

\noindent\textbf{Proof of Theorem \ref{theo095815052015}.} This time, we shall use Lemma \ref{gaussgen1203060502015hd} to show that for any $\bA_1,\ld,\bA_r$, $N\ti N$ deterministic matrix of the form
$$
\bA_m \ = \ \bpm \bA_{m,p} & 0 \\ 0 & 0 \epm, \ \quad \bA_{m,p} = \big( a^m_{i,j}\big)_{i,j=1}^p \in \M_p(\C),
$$
the vector 
$$
\lf( \ff N \tr \big(\bU^* \bT_1 \bU \bA_1 \big) ,\ld, \ff N \tr \big(\bU^* \bT_r \bU \bA_r \big)\ri)
$$
converges weakly to a Gaussian multivariate. Thanks to Theorem \ref{lem095815052015}, we know that for each $m$, 
$$
\ff N \tr\big(\bU^* \bT_m \bU \bA_m\big)
$$
is asymptotically Gaussian. Then, we show that 
$$
\Ec{\lf(\ff N \tr\big(\bU^* \bT_1 \bU \bA_1\big)\ri)^{p_1}\lf( \ff N \tr\big(\bU^* \bT^*_1 \bU \bA^*_1\big)\ri)^{q_1} \cdots \lf(\ff N \tr\big(\bU^* \bT_r \bU \bA_r\big)\ri)^{p_r}\lf( \ff N \tr\big(\bU^* \bT^*_r \bU \bA^*_r\big)\ri)^{q_r}}
$$
satisfies \eqref{cond105062015} and \eqref{cond205062015} using Lemma \ref{lem151727052015}.
\hfill$\square$\\

\noindent\textbf{Proof of Lemma \ref{lem151727052015}.} We know from \cite[Proposition 3.4]{MingoSniady}
\beq
\Ec{\prod_{t=1}^q  \sqrt N \tr\big( \bU^* \bT_t \bU \bA_t\big)}&=&N^{q/2} \sum_{\sigma,\tau \in S_q} \wg(\tau\circ\sigma^{-1})\tr_\tau\big(\bA_i\big)_{i=1}^q \tr_\sigma\big( \bT_i\big)_{i=1}^q
\eeq
where $\wg$ is a function called  the \emph{Weingarten function}. Moreover, for $\sigma \in S_q$,   the asymptotical behavior of $\wg(\sigma)$ is at most given by
\begin{eqnarray} \label{wg20203047840203020106}
\wg(\sigma) &=& O\left( N^{-q}\right).
\end{eqnarray}

First, one should notice that if $\sigma$ has one invariant point (which means a cycle of size one in its cycle decomposition), then
$$
\tr_\sigma\big( \bT_i\big)_{i=1}^q \ = \ 0,
$$
also, if $\sigma$ has $r$ cycles in its cycle decomposition, then, by the Holder inequality,
$$
\tr_\sigma\big( \bT_i\big)_{i=1}^q \ = \ \OO{N^r}.
$$
Actually, the maximum of cycles in its decomposition that can have $\sigma$ without any $1$-sized cycle is $\ds\lf\lfloor\f{q}{2} \ri\rfloor$ so that, using \eqref{wg20203047840203020106}
\beq
N^{q/2}\wg ( \sigma \circ \tau^{-1})\tr_{\tau}\big( \bA_i\big)_{i=1}^q \tr_\sigma \big( \bT_i\big)_{i=1}^q & = & \OO{N^{ \lf\lfloor\f{q}{2} \ri\rfloor-\f{q}{2}}},
\eeq
so that first, if $q$ is odd
\beq
\Ec{\prod_{t=1}^q  \sqrt N \tr\big( \bU^* \bT_t \bU \bA_t\big)} & = & \oo1.
\eeq
Moreover, if $q = 2r$, then the only way to have 
$$
N^{q/2}\wg ( \sigma \circ \tau^{-1})\tr_{\tau}\big( \bA_i\big)_{i=1}^q \tr_\sigma \big( \bT_i\big)_{i=1}^q \ \neq \ \oo1
$$
is to have \bgt
\ite $\tau = \sigma$, 
\ite $\sigma$ is a product of $\ds\f{q}{2}=r$ transpositions with disjoint support.
\ent
One easily conclude. \hfill$\square$

\subsection{Moments of a complex Gaussian variable.} The following lemma allows to prove that a random variable is Gaussian if and only if its moments satisfy an induction relation.

\begin{lemme}\la{lemgauss15182015}
Let $Z$ be a complex Gaussian variable such that
\beqy \la{cond11444201505}
\Ec{Z} \ = \ 0, \ \ \ \Ec{Z^2} \ = \ \tau^2, \ \ \ \Ec{|Z|^2} \ = \ \sigma^2.
\eeqy
Then, for all $p \geq 1$ 
\beqy\la{cond21444201505}
\Ec{Z^{2p}} \ = \ p!! \tau^{2p} & \text{ and } & \Ec{Z^{2p+1}} \ = \ 0, \qquad \qquad \big(\text{where }  p!! := \f{(2p)!}{2^p p!} \big)
\eeqy
 also, for all $p,q \geq 0$, 
\beqy\la{cond31444201505}
\Ec{Z^{p+2} \ol{Z}^{q+2}} & = & \sigma^2 (q+2) \Ec{Z^{p+1} \ol{Z}^{q+1}} + \tau^2 (p+1)\Ec{Z^{p} \ol{Z}^{q+2}} \\ & = & \sigma^2 (p+2) \Ec{Z^{p+1} \ol{Z}^{q+1}} + \ol\tau^2 (q+1)\Ec{Z^{p+2} \ol{Z}^{q}}. \nonumber
\eeqy
Conversely, any complex random variable $Z$ satisfying \eqref{cond11444201505},\eqref{cond21444201505} and \eqref{cond31444201505} is a complex Gaussian variable.
\end{lemme}
\bpr
First, recall that if $Z=X_1+\ii X_2$ is a complex random Gaussian such that
$$
\Ec{Z} \ = \ 0, \ \ \ \Ec{Z^2} \ = \ \tau^2, \ \ \ \Ec{|Z|^2} \ = \ \sigma^2,
$$
then, its Fourier transform is given, for $t=t_1+\ii t_2 \in \C$, by
\beq
\Phi_Z(t) & :=& \E \exp\lf(i(X_1 t_1 + X_2 t_2) \ri) \\
& = &  \exp\lf(-\ff 4 \big( (t_1^2 + t_2^2)\sigma^2 + (t_1^2-t_2^2)\re(\tau^2) + 2t_1 t_2 \im(\tau^2)\big)\ri)
\eeq
We define the differential operators
\beqy
\partial_t \ : = \ \partial_1 + \ii \partial_2 & ; & \partial_{\ol t} \ : = \ \partial_1 - \ii \partial_2
\eeqy
so that
\beqy
\Ec{Z^p \ol{Z}^q} & = & (-\ii)^{p+q} \partial_t^p \partial_{\ol t}^q \Phi(t)\big|_{t=0} .
\eeqy
One can easily compute

\beq
&& \partial_t \Phi(t) \ = \ - \ff 2 \lf( t \sigma^2 + \ol t \tau^2\ri) \Phi(t) \quad ; \quad \partial_{\ol t} \Phi(t) \ = \ - \ff 2 \lf( \ol t \sigma^2 +  t \ol\tau^2\ri) \Phi(t) \\
\eeq
therefore, for any $p \geq 0$, $q \geq 0$,
\beq
\partial_t^{p+2} \Phi(t) & = & \partial^{p+1}_t \lf(- \ff 2 \lf( t \sigma^2 + \ol t \tau^2\ri) \Phi(t)\ri) \\
& = & -\ff 2 t \sigma^2 \partial^{p+1}_t \Phi(t) - \ff 2 \tau^2 \partial^{p+1}_t \big(\ol t \Phi(t) \big)\\
& = &  -\ff 2 t \sigma^2 \partial^{p+1}_t \Phi(t) - \ff 2 \tau^2 \ol t \partial^{p+1}_t \Phi(t) - \tau^2 (p+1)\partial_t^{p}\Phi(t)
\eeq
and
\beq
\partial_t^{p+2}\partial_{\ol t}^{q+2} \Phi(t) & = & \partial_t^{p+1}\partial_{\ol t}^{q+2} \lf(- \ff 2 \lf( t \sigma^2 + \ol t \tau^2\ri) \Phi(t)\ri) \\
& = &  -\ff 2 \sigma^2 \partial_{\ol t}^{q+2} \big( t \partial_t^{p+1} \Phi(t)\big) - \ff 2 \tau^2 \partial_{t}^{p+1} \big( \ol t \partial_{\ol t}^{q+2} \Phi(t)\big)\\
& = & -  \sigma^2 \lf( \f{t}{2} \partial_t^{p+1} \partial^{q+2}_{\ol t} \Phi(t) + (q+2) \partial_t^{p+1} \partial^{q+1}_{\ol t} \Phi(t)\ri) \\ && \ \quad - \ \tau^2 \lf(\f{\ol t}{2}\partial_t^{p+1} \partial^{q+2}_{\ol t} \Phi(t) +(p+1) \partial_t^{p} \partial^{q+2}_{\ol t} \Phi(t) \ri),
\eeq
hence,
\beq
\Ec{Z^{p+2}} & = & \tau^2 (p+1) \Ec{Z^p}, \\
\Ec{Z^{p+2} \ol{Z}^{q+2}} & = & \sigma^2 (q+2) \Ec{Z^{p+1} \ol{Z}^{q+1}} + \tau^2 (p+1)\Ec{Z^{p} \ol{Z}^{q+2}}
\eeq
and the same way, 
\beq
\Ec{\ol Z^{p+2}} & = & \ol\tau^2 (p+1) \Ec{\ol Z^p},\\
\Ec{Z^{p+2} \ol{Z}^{q+2}} & = & \sigma^2 (p+2) \Ec{Z^{p+1} \ol{Z}^{q+1}} + \ol\tau^2 (q+1)\Ec{Z^{p+2} \ol{Z}^{q}}.
\eeq
Conversely, one can easily prove by induction that any complex random variable $Z$ satisfying \eqref{cond11444201505},\eqref{cond21444201505} and \eqref{cond31444201505} has all its moments uniquely determined and since the complex Gaussian variable also satisfies \eqref{cond11444201505},\eqref{cond21444201505} and \eqref{cond31444201505}, one can conclude.  

\epr

More generally, one can show the following lemma

\begin{lemme} \la{gaussgen1203060502015hd}
Let $(X_1,\ld,X_r)$ be a centered complex Gaussian vector. Then, for all non negative integers $p_1, q_1 ,\ld,  p_r,q_r$, for all $i \in \{1,\ld,r\}$ 
\beqy\label{cond105062015}
\text{if } p_i \geq  1, \qquad \Ec{X_1^{p_1} \ol{X_1}^{q_1}\cdots X_r^{p_r}\ol{X_r}^{q_r}} & = & (p_i-1)\Ec{X_i^2} \Ec{X_1^{p_1} \ol{X_1}^{q_1}\cdots X_i^{p_i-2} \ol{X_i}^{q_i}\cdots X_r^{p_r}\ol{X_r}^{q_r}} \nonumber\\
& + & \sum_{\ds^{j=1}_{j \neq i }}^r p_j \Ec{X_i X_j} \Ec{X_1^{p_1} \ol{X_1}^{q_1}\cdots X_i^{p_i-1}\cdots X_j^{p_j-1}\cdots X_r^{p_r}\ol{X_r}^{q_r}} \\
& + &  \sum_{j=1}^r q_j \Ec{X_i \ol{X_j}} \Ec{X_1^{p_1} \ol{X_1}^{q_1}\cdots X_i^{p_i-1}\cdots \ol{X_j^{q_j-1}}\cdots X_r^{p_r}\ol{X_r}^{q_r}} \nonumber
\eeqy
\beqy\label{cond205062015}
\text{if } q_i \geq  1, \qquad  \Ec{X_1^{p_1} \ol{X_1}^{q_1}\cdots X_r^{p_r}\ol{X_r}^{q_r}} & = & (q_i-1)\Ec{\ol{X_i^2}} \Ec{X_1^{p_1} \ol{X_1}^{q_1}\cdots X_i^{p_i} \ol{X_i}^{q_i-2}\cdots X_r^{p_r}\ol{X_r}^{q_r}} \nonumber\\
& + & \sum_{{j=1}}^r p_j \Ec{\ol{X_i} X_j} \Ec{X_1^{p_1} \ol{X_1}^{q_1}\cdots \ol{X_i^{q_i-1}}\cdots X_j^{p_j-1}\cdots X_r^{p_r}\ol{X_r}^{q_r}}\\
& + &  \sum_{\ds^{j=1}_{j \neq i}}^r q_j \Ec{\ol{X_i} \ol{X_j}} \Ec{X_1^{p_1} \ol{X_1}^{q_1}\cdots \ol{X_i^{q_i-1}}\cdots \ol{X_j^{q_j-1}}\cdots X_r^{p_r}\ol{X_r}^{q_r}}\nonumber
\eeqy
with the convention that $X^{-1} = 0$.\\
Conversely, if $X_1,\ld,X_r$ are $r$ centered Gaussian variables satisfying \eqref{cond105062015} and \eqref{cond205062015}, then $(X_1,\ld,X_r)$ is a centered complex Gaussian vector.
\end{lemme}

\bpr
In the same spirit as the proof of Lemma \ref{lemgauss15182015}, we obtain \eqref{cond105062015} and \eqref{cond205062015} by derivating the Fourier transform. The converse is proved  by induction.
\epr

\end{document}